\documentclass[]{article}
\usepackage{amsmath}
 \usepackage{amssymb}
\usepackage{graphicx}
\def\bel{\begin{equation}\label}
\def\eeq{\end{equation}}
\newtheorem{Definition}{Definition}[section]
\newtheorem{Theorem}{Theorem}[section]
\newtheorem{Remark}{Remark}[section]

\newtheorem{Proposition}{Proposition}[section]
\newtheorem{Corollary}{Corollary}[section]
\newtheorem{Example}{Example}[section]

\title{Asymptotic problems in optimal control  with a vanishing Lagrangian and  unbounded data}
\author{M.Motta,  C.Sartori\\
Dipartimento di Matematica  \\
Via Trieste, 63 - 35121 Padova, Italy\\
Telefax (39)(049) 8271428\\
e-mail: motta@math.unipd.it \\ sartori@math.unipd.it}
\def\fudge{\mathchoice{}{}{\mkern.5mu}{\mkern.8mu}}
\def\bbc#1#2{{\rm \mkern#2mu\vbar\mkern-#2mu#1}}
\def\bbb#1{{\rm I\mkern-3.5mu #1}}
\def\bba#1#2{{\rm #1\mkern-#2mu\fudge #1}}
\def\bb#1{{\count4=`#1 \advance\count4by-64 \ifcase\count4\or\bba A{11.5}\or
\bbb B\or\bbc C{5}\or\bbb D\or\bbb E\or\bbb F \or\bbc G{5}\or\bbb H\or
\bbb I\or\bbc J{3}\or\bbb K\or\bbb L \or\bbb M\or\bbb N\or\bbc O{5} \or
\bbb P\or\bbc Q{5}\or\bbb R\or\bbc S{4.2}\or\bba T{10.5}\or\bbc U{5}\or%
\bbb P\or\bbc Q{5}\or\bbb R\or\bba S{8}\or\bba T{10.5}\or\bbc U{5}\or
\bba V{12}\or\bba W{16.5}\or\bba X{11}\or\bba Y{11.7}\or\bba Z{7.5}\fi}}
\def \ol{\overline}
\def\qed{$\square$}
\def \t{\tau}

\def \R {{\bb R}}
\def \C {{\mathcal  T}}
\def \N {{\bb N}}

\def \vv{\vskip 0.3 truecm}

\def \vsm{\vskip 0.2 truecm}
\begin{large}
\begin{document}
\maketitle
\begin{abstract}{In this paper we give a representation formula for the limit of the finite horizon problem as the horizon becomes infinite,  with a nonnegative Lagrangian and unbounded data. It is related to the limit of the  discounted infinite horizon  problem,  as the discount factor goes to zero.  We give sufficient conditions to characterize the limit function as unique nonnegative solution of the associated HJB equation. We also briefly discuss the ergodic problem.  }
\end{abstract}

\footnotetext  {$({\bf *})$ This research is partially supported by the Marie Curie ITN SADCO, FP7-PEOPLE-2010-ITN n. 264735-SADCO and by the Gruppo Nazionale per lÕAnalisi Matematica, la Probabilit\`a e le loro Applicazioni (GNAMPA) of the Istituto Nazionale di Alta Matematica (INdAM). }

\section{Introduction}

The main goal of this paper is to discuss, in the case of a vanishing Lagrangian $l\ge0$ and  truly unbounded data and controls,  the limit as $t$ tends to $+\infty$  of
{\it  the finite horizon value function}
$$
{\mathcal V}(t, x)\doteq\inf_{\alpha}\int_0^tl(y(\t),\alpha(\t))\,d\t,
$$
   and the limit  as $\delta$ tends to $0^+$   of
{\it  the discounted infinite horizon value function}
  $$
{\mathcal V_\delta}( x)\doteq\inf_{\alpha}\int_0^{+\infty}e^{-\delta\,t}l(y(\t),\alpha(\t))\,d\t,
$$
where $f$, $l$ are given functions,   $\alpha(\t)\in A\subset\R^m$ is the control and the trajectory is given by $\dot y(\tau)= f(y(\tau),\alpha(\tau))$, $y(0)=x$.  

\vv
These limits have been extensively studied in the literature.  On the one hand, the approximability of the {\it infinite horizon value function}
$$
{\mathcal V}( x)\doteq\inf_{\alpha}\int_0^{+\infty}l(y(\t),\alpha(\t))\,d\t,
$$
by the  finite horizon value functions is classically required in most applications (see \cite{CHL}) and it also represents  the key point of several comparison results by viscosity solution methods. On the other hand, recently a lot of work has been devoted to  the study of  the two  ergodic limits  $\lim_{t\to+\infty} {\mathcal V}(t,x)/t$ and $\lim_{\delta\to 0^+}\delta\, {\mathcal V}_{\delta}(x)$.
  We refer to \cite{BCD} for a presentation of the basic results in the deterministic case, and to \cite{AL} for the stochastic case. The same questions have been  addressed  in $L^\infty$ control problems (see  \cite{AB} and the references therein). 

The main novelty of this paper  is the generality of the hypotheses under which the results are obtained, suitable to a wide range of applications in the framework of optimal control theory.  Precise assumptions will be stated in Section \ref{Ass}, here we just point out that  we can consider coercive and non coercive nonnegative  Lagrangians, with arbitrary growth in the state variable and without restrictions on the set 
$$
{\mathcal Z}\doteq\{x: \ l(x,a)=0 \ \text{for some $a$}\}.
$$
 For instance, the dynamics 
can be control-affine,  $f(x,a)=f_0(x)+\langle F(x),a\rangle$,  where $f_0$, $F$ are locally Lipschitz functions with linear growth in $x$.
In particular we cover    (nonlinear generalizations of) LQR  problems with  $l(x,a)=x^TQx+a^TRa$, where $Q$ and $R$ are  symmetric matrices,   $R$ is positive definite and $Q$ is  positive semidefinite. We can also  allow for control-affine Lagrangians,  $l(x,a)=l_0(x)+l_1(x)|a|$ with $l_0\ge0$, $l_1>0$ continuous and with arbitrary growth in $x$, used in some economics models, mostly in singular stochastic control (see \cite{FS} and the references therein).  
 
We  show that the function  $\Sigma(x)\doteq\underset{t\to+\infty}\lim{\mathcal V}(t,x)$ is l.s.c. and we characterize it   as the minimal nonnegative supersolution to the limit HJB equation at every $x$ where it is finite. The representation formula,  when $A$ is compact,  is given, as expected, by the value function of the so-called {\it relaxed} infinite horizon problem.  Adding some mild assumptions on the data,  it is also equal to the l.s.c. envelope of the infinite horizon value function, ${\mathcal V}_*(x)$.
 
When  $A$ is unbounded,  the relaxed problem is not defined. In this case, we can still give a representation formula for $\Sigma$ by  introducing  an {\it extended} infinite horizon problem,   which has a compact control set.  Denoting by $V$ the value function of the  extended problem, we prove that  $\Sigma$ coincides with the relaxed version of $V$ and also with its l.s.c. envelope, $V_*$, under the same assumptions as for $A$ compact.  
  In particular, in classical impulsive control problems,  the extended setting is equivalent to the replacement of controls with measures. In Theorem \ref{extT} we give sufficient conditions to have $V$ equal to ${\mathcal V}$.
  
  \noindent  We obtain the same characterizations for $\underset{\delta\to0^+}\lim{\mathcal V}_\delta(x)$, assuming   ${\mathcal V}_\delta$   bounded. 
  
  In general,   $\Sigma$ is  not u.s.c. and the limit HJB equation does not have a unique solution. We give explicit sufficient conditions under which   $\Sigma$ turns out to be continuous and the unique nonnegative solution to the  HJB equation.
\vsm

We spend a few words on the ergodic problem.  Starting from the papers   \cite {AL} and \cite{A}, a huge amount of literature has been devoted to the subject,   initially in the case of bounded domains or periodic data and under some global controllability assumptions. The first results  have been developed and generalized in several directions (see e.g.  \cite {BR}, \cite{GLM}, \cite{QR},   and the references therein).
 Here we focus our attention mainly on the case where the set ${\cal Z}\ne\emptyset$ and the infinite horizon value function is finite,  case in which the ergodic limits turn out to be zero.   We limit ourselves to  showing  how it is possible, under periodicity of the data and a complete controllability condition, to obtain the results of \cite{A} in our framework.

Some final bibliographical remarks.  In this paper, we extend  to the dynamics and Lagrangians described above,  many results already proved when some of the data of the problem are bounded. In so doing we get some results new also for the compact control case.  When the control set is unbounded, our approach is based on a compactification method introduced in 
\cite{BrRa} (see also \cite{MiRu});  for a more complete survey we refer to  \cite{BP} and  the references therein.
In particular, the finite horizon problem with  both coercive and weakly coercive Lagrangians was treated in \cite {RS}, while  exit-time problems with a nonnegative Lagrangian were investigated  in \cite {MS}. Moreover some optimality principles were extended in \cite{M}  to the HJB equations involved in several  optimal control problems of this kind. This approach has also been applied  to some stochastic control problems (see e.g. \cite{MS2} and the references therein).
 
 
In Section \ref{Ass} we state the problem precisely. In Section  \ref{Gen}  we introduce the extended setting for $A$ unbounded and give sufficient conditions in order to have the extended infinite horizon value function coinciding with ${\mathcal V}( x)$; then we define the  relaxed and the relaxed extended problems. 
Section \ref{FHA} is devoted to characterize the limit as $t$ tends to $+\infty$ of the finite horizon value functions, while the limit as $\delta$ tends to $0^+$ of the discounted value functions is studied in Section \ref{Dis}.  In Section \ref{Un} we state a uniqueness result for  the solution of the limit HJB equation. The ergodic problem is investigated in Section \ref{erg}.  The discounted and the ergodic problems have been treated in the last two sections, since they  are studied under assumptions  not required for the previous results.  
 
 \noindent{\bf Notations.} 
For any function 
$u:\ R^n\to\ R\cup\{+\infty\}$, we will denote the set $\{x\in\ R^n: \ u(x)<+\infty\}$ by $Dom(u)$. $\R_+\doteq[0,+\infty[$.
A function $\omega:\R_+\times\R_+\to\R_+$
is called a {\it modulus} if: $\omega(\cdot,R)$ is increasing in a neighborhood of $0$, 
continuous at $0$, and $\omega(0, R)=0$ for every $R>0$; $\omega(r,\cdot)$ is increasing for every $r$.  Let $D\subset \ R^N$ for some $N\in\N$.  $\forall r>0$  we 
will denote by $D_r$ the closed set $\overline{B(D,r)}$, while $D_r^c=\R^N\setminus D_r$. Moreover, $\chi_D$ will denote the characteristic function
of $D$, namely for any $x\in \R^N$ we set $\chi_D(x)=1$ if $x\in D$ and $\chi_D(x)=0$ if $x\notin D$.  

%
 
\section{  Assumptions and statement of the  problem }\label{Ass}
We consider  a  nonlinear control  system having the form 
\begin{equation}\label{S}
\dot y(\tau)= f(y(\tau),\alpha(\tau)), \qquad y(0)=x  
\end{equation}
 and an undiscounted  payoff 
\begin{equation}\label{1.2}
{\mathcal J}(t,x,\alpha)=\int_0^tl(y(\t),\alpha(\t))\,d\t, 
\end{equation}
where $\alpha(\t)\in A\subset\R^m$,    and $l$ is nonnegative.   
For any $x\in\R^n$, we define the infinite horizon value function 
\begin{equation}\label{infty}
 {\cal V}(x)\doteq \inf_{\alpha\in{\mathcal{ A}}}{\mathcal J}(+\infty, {x},\alpha ), 
\end{equation}
where the admissible controls set ${\mathcal{ A}}$ is given by (\ref{camm}) below.
\vv
The following  hypotheses  (H0), (H1)   will be assumed throughout the whole paper. 
{\it \begin{itemize}
\item[{\bf (H0)}]   The control set $A\subset\ R^m$  is either compact or  a convex, closed, nontrivial  cone containing the origin. 
\item[] The functions $f:\R^n\times A\to\R^n$, $l:\R^n\times A\to\R$  are continuous; there exist $p$, $q\in\N$, \  $q\ge p\ge1$, $M>0$, and for any $R>0$ there are $L_R$, $M_R>0$ and a modulus   $\omega(\cdot,R)$, such that $\forall x$, $x_1$,  $x_2\in\R^n$, $\forall a\in A$,
\begin{equation}\label{hp0}
\aligned
|f(x_1,a)-f(x_2,a)|&\le L_R(1+|a|^p)|x_1-x_2|,\\
 |l(x_1,a)-l(x_2,a)|&\le (1+|a|^q)\,\omega(|x_1-x_2|,R) \\
  0&\le l(x,a)\le M_R(1+|a|^q) \quad \text{if }  |x_1|, \, |x_2|,\,  |x|\le R, \\
 |f(x,a)|&\le M(1+|a|^p)(1+|x|).
\endaligned
\end{equation}
\end{itemize}}   

\noindent If $A$ is compact, the above assumptions  reduce to the continuity of $l$  and to  the usual hypotheses of sublinear growth and local Lipschitz continuity in $x$, uniformly w.r.t. $a$,  for $f$. With a small abuse of notation, in this case we will denote again by $L_R$ the quantity  $\max\{L_R(1+|a|^p): \ a\in A\}$ and similarly for the other constants appearing in (H0). 
 
 \vv
 When $A$ is unbounded, we will always assume  at least weak coercivity  together with a regularity hypothesis in the control variable at infinity:

 {\it \begin{itemize}
\item[{\bf (H1)}]
 There exist some constants $C_1\ge0$, $C_2>0$ such that
\begin{equation}\label{CC}
l(x,a)\ge C_2|a|^q-C_1  \qquad \forall (x,a)\in\R^n\times A
\end{equation}
  and  $q\ge p$, where $q$ and $p$ are the same as in {\rm (H0)}.  
 
  Let $\Phi\in\{f,l\}$. There  exists   a continuous function $\Phi^{\infty}$, called  the 
{\rm recession function of $\Phi$}, verifying 
\bel{rec}
\lim_{\rho\to 0^+}\rho^{q}\Phi(x,\rho^{-1}a)\doteq \Phi^{\infty}(x,a)
\eeq
 {\rm uniformly} on compact sets of $\ R^n\times A$. 
\end{itemize} }

  Condition  (\ref{CC}),  for $q>p$  is known as {\it coercivity} and it is used to yield suitable compactness properties for the set of the admissible controls. It is satisfied,  for instance,  in the LQR problems anticipated in the Introduction.   If $q=p$, instead, (\ref{CC}) is sometimes called {\it weak coercivity.} In this case  the natural framework of all our optimization problems is that of generalized or {\it impulsive} controls, since minimizing sequences of trajectories may converge to a discontinuous function. In Section \ref{Gen} the generalized setting will be introduced  in terms of some extended problems.  This approach is suitable to study, for instance,  problems in which both the dynamics and the Lagrangian are control-affine.
   
  \begin{Example}{\rm Functions  $f$ and $l$ which are  polynomials  in the control variable $a$, admit the recession function introduced in (\ref{rec}).  If, for instance,   there are some continuous functions  $f_i$,  $F_{ij}$  such that  
  $$ 
  f(x,a)= f_0(x)+ \sum_{i=1}^m f_i(x)a_i+\sum_{i,\,j=1}^m F_{ij}(x)a_i\,a_j \qquad \forall (x,a)\in\R^n\times A,
  $$
 $p=2$ and  $f^\infty(x,a)=\sum_{i,\,j=1}^m F_{ij}(x)a_i\,a_j$   if $q=2$;  $f^\infty(x,a)\equiv 0$ if $q>2$.
  }
  \end{Example}
  
  \noindent Notice that if $q>p$, then one always has  $f^\infty\equiv0$.

 \vv 
Let ${\mathcal B}$ denote the set of the
 Borel--measurable functions. The controls  $\alpha$ are assumed to belong to the set      
\bel{camm}
{\mathcal A}\doteq {\mathcal{ B}}\cap L^q_{loc}(\R_+,A),
\eeq
coinciding with $\mathcal{ B}$ when $A$ is compact.   
  For any $x\in \R^n$  and   for any  control $\alpha\in{\mathcal A}$,   (\ref{S}) admits just one solution,  defined on the whole interval $\R_+$. We use  $y_x(\cdot,\alpha)$   to  denote such a  solution.  When $A$ is unbounded the  control set ${\mathcal A}$ is the largest set where both payoff and trajectory are surely defined for all $t\ge0$.   In fact, in view of the  coercivity condition (\ref{CC}) (weak, if $q=p$),  such a choice is not a restriction, since for any   measurable control $\alpha$,
$$
J(t,x,\alpha)\ge C_2\int_0^t|\alpha(\t)|^q\,d\t-C_1 t \qquad \forall t>0
$$
so that for controls $\alpha\notin{\mathcal A}$ we will never obtain a  finite cost. 
In particular,  if $C_1=0$  we can consider merely controls in $L^q(\R_+,A)$.

\vv
  Let us write two estimates, useful in the sequel, that can be obtained by standard tools.   
For every $x, z\in\R^n$, $\forall\alpha\in{\mathcal A}$, and   $\forall t\ge0$ one has 
 \bel{est1}|y_x(t,\alpha)|\le\left(|x|+Mt+M\int_0^t|\alpha(t')|^p\,dt'\right){\rm e}^{M (t+\int_0^t|\alpha(t')|^p\,dt')}\eeq
and, if $\exists R>0$ such that $|y_x(t',\alpha)|, |y_z(t',\alpha)|\le R$\, $\forall t'\in[0,t]$, then
 \bel{est2}
 |y_x(t,\alpha)-y_z(t,\alpha)|\le|x-z|{\rm e}^{L_R(t+\int_0^t|\alpha(t')|^p\,dt')}.
 \eeq

 \vv
For some results we will use  the following hypothesis (H2). 
  {\it \begin{itemize}
\item[{\bf (H2)}]   There is some nonempty closed set  $\C\subset\R^n$ with compact boundary such that 
${\mathcal V}(x)=0$ for any $x\in\C$ and
\bel{contF}
 \lim_{x\to\bar x}{\mathcal V}(x)=0 \qquad \forall\bar x\in\partial\C.
\eeq
\end{itemize}}

\begin{Remark}{\rm  
Assume that  ${\mathcal V}(x)\le \int_0^{+\infty}l(y(\t),\alpha(\t))\,d\t<+\infty$ for some $x$ and a control $\alpha\in{\mathcal A}$   with $|y_x(t)|\le \bar R$ for all $t\ge0$ and some $\bar R>0$. Then it is not difficult to show that there exists $\bar x$ where ${\mathcal V}_*(\bar x)=0$.  Therefore, if ${\mathcal V}$ is continuous at $\bar x$ then $\bar x\in{\mathcal T}$  and hypothesis (H2) holds at $\bar x$.

The hypothesis ${\mathcal V}\equiv 0$ in $\C$ is satisfied, e.g.,  if   $\C\times\{0\}$ is a viability set for the vector field $(f,l)$.\footnote{Let $F(x)\doteq \{(f(x,a),l(x,a)): \ a\in    A\}$.     Any closed subset $K\subset\R^n\times\R$ will be called   a {\it viability set} for $(f,l)$ if  for any $(x_0,\lambda_0)\in K$ there is a solution $(y,\lambda)$ of 
the differential inclusion
$$
(\dot y(t), \dot\lambda(t))\in F(y(t))\qquad   t\ge 0
$$
 such that  $(y(0),\lambda(0))=(x_0,\lambda_0)$ and $(y(t),\lambda(t))\in K$ \ $\forall t>0$ (see   \cite{AF}).} Sufficient viability conditions can be found e.g. in \cite{AF}.}
 \end{Remark}
 
   As shown in \cite{MR}, a sufficient condition for  (\ref{contF}) is the existence of a local MRF $U$, defined, for the case $A$ compact, as follows.
 \begin{Definition}\label{CMRL}{\rm [MR]}   Given an open set $\Omega\subset\R^n$, $\Omega\supset\C$  we say that $U:\Omega\setminus\overset{\circ} \C \to\R_+$ is a {\rm local Minimum Restraint Function,}  in short, a local  MRF  for $l$,  if $U$  is  continuous on $\Omega\setminus\overset{\circ}\C$, locally  semiconcave, positive definite,  proper\,\footnote{$U$ is said {\it positive definite on $\Omega\setminus\C$} if  $U(x)>0$ \,$\forall x\in\Omega\setminus\C$ and $U(x)=0$ \,$\forall x\in\partial\C$.  $U$ is called {\it proper on $\Omega\setminus\C$} if $U^{-1}(K)$ is compact for every compact set  $K\subset\R_+$.}  on $\Omega\setminus\C$, $\exists U_0\in]0,+\infty]$ such that
 $$
\lim_{x\to x_0, \  x\in\Omega}U(x)=U_0 \ \  \forall x_0\in\partial\Omega;  \quad
U(x)<U_0 \quad\forall x\in\Omega\setminus\overset{\circ} \C,
$$
and, moreover,   $\exists{{k}}> 0$ such that, for every $x\in \Omega\setminus\C$,  
\begin{equation}\label{C1}
\min_{ a\in A} \left\{   \langle p, f(x,a)\rangle +k\,l(x,a) \right\}<0
 \qquad \forall p\in D^* U(x),
\end{equation}
   where  $D^* U(x)$ is   the    set of limiting gradients of $U$ at $x$. 
\end{Definition}

 For the case $A$ unbounded, as proved in  Remark 2.5 of   \cite{MS},
 we can  consider  the following hypothesis: 
\vv
\noindent {\it  There exists a local MRF $U$ for $l$ such that $\forall x\in \Omega\setminus\C$:
\begin{equation}\label{C1unb}
\min_{ a\in A\cap \overline{B(0,R(U(x)))}} \left\{   \langle p, f(x,a)\rangle +k\,l(x,a) \right\}<0,
 \qquad \forall p\in D^* U(x),
\end{equation}
where $R:\,\,]0,\sigma]\to]0,+\infty[$ is a decreasing continuous function (in particular,  we may have $\lim_{\delta\to0^+}R(\delta)=+\infty$).}

Let us observe that any MRF is a Control Lyapunov function for the system w.r.t. $\C$, which yields local asymptotic controllability to $\C$. For the notions borrowed from nonsmooth analysis,  we refer to  \cite{CS}.

\section{Generalized and relaxed control problems} \label{Gen}

Following  the so called graph-completion approach proposed in \cite{BrRa}, as developed in \cite{RS},  when $A$ is unbounded   we  represent   generalized controls  and  trajectories  as   reparametrizations (through a   time-change, possibly discontinuous in case $q=p$)  of controls and trajectories of  the   extended minimization problems below, involving bounded-valued controls.   
Then we investigate the well-posedness of the generalized setting, that is, when  the infima  over ordinary and generalized controls are the same. We do this for both  the finite and for the infinite horizon problem.  Let us remark that dealing with a compact set of controls as the generalized control set is, has two main advantages. On the one hand,  it allows to introduce the relaxed problem for which an optimal control   exists.
On the other hand,  the relative Hamiltonian, differently from the original, is continuous and satisfies some crucial growth and regularity properties.
The exploitation of both these aspects yields many results.
  
\subsection{Generalized problems and well posedness}
 Throughout this subsection we assume  $A$  unbounded.  Let us define  on $\R^n\times (\R_+\times A)$ the {\it extended} dynamics and Lagrangian  $\overline{f}$,  $\overline{l}$  as follows:
\begin{equation}\label{extended data}
\overline{\Phi}(x,w_0,w)\doteq\left\{\begin{array}{l}
w_0^q \, \Phi(x,w_0^{-1}w) \ \ {\rm if \ }w_0\ne0 \\ \Phi^\infty(x,w) \quad {\rm if \ }w_0=0.
\end{array}\right.   \qquad \Phi\in\{f,l\},
\end{equation}
where $\Phi^\infty$ is defined in (H1).  $\overline{f}$,  $\overline{l}$ are continuous, $q$-positively homogeneous in the control variable $(w_0,w)$ and inherit    properties analogous  to those of  $f$ and $l$, respectively (see e.g. \cite{M}). 

 Let $S(A)\doteq (\R_+\times A)\cap \{(w_0,w): \ w_0^q+|w|^q= 1\}$.  Define the set of {\it extended controls}  as
\begin{equation}\label{Gamma}
\Gamma\doteq\left\{(w_0,w): \  \  (w_0,w)\in {\mathcal{ B}}(\R_+,S(A)) \right\},
\end{equation} 
and    $\forall (w_0,w)\in\Gamma$  denote by $\xi(\cdot)\equiv \xi_{ x}(\cdot, w_0,w)$ the {\it extended trajectory} solving the {\it extended control system}
\begin{equation}\label{SE}
 \xi'(s)=\overline f(\xi(s),w_0(s),w(s)) \qquad \xi(0)=x.  
 \end{equation}
For any $S>0,$ the {\it extended payoff} is given by
   \begin{equation}\label{Je}
J(S,x ,w_0,w)=\int_0^{S}\overline l(\xi(s),w_0(s),w(s))\,ds.
\end{equation} 

As recalled in Proposition \ref{REmb} below, the solutions to (\ref{SE}) are simply  time-reparametrizations of trajectories of (\ref{S}) if the controls belong to 
\bel{G+}
\Gamma^+\doteq\Gamma\cap\left\{ (w_0,w): \  w_0>0 \text{ a.e.}\right\}.
\eeq
 \begin{Proposition}\label{REmb}  {\rm [MS]}
   For any $\alpha\in{\mathcal A}$  let us define
   $s(t)\doteq\int_0^t(1+|\alpha(\tau)|^q)\,d\tau$ for all $t\ge0$ and denote by $t:\R_+\to [0,+\infty[$ its inverse function. Then   $(w_0,w)$ defined by  
$w(\cdot)\doteq\frac{\alpha(t(\cdot))}{{(1+|\alpha(t(\cdot)|^q})^{{1}/{q}}}$, \, $w_0(\cdot)\doteq (1-|w(\cdot)|^q)^{1/q}$, 
   belongs to $\Gamma^+$ and $y_x(t(\cdot),\alpha)$ is the solution of (\ref{SE})  associated to $(w_0,w)$.
   
\noindent  Vice-versa, for any 
 $(w_0,w)\in\Gamma^+$ such that 
 \bel{w-0}
\int_0^{+\infty} w_0^q(s)\,ds=+\infty,
\eeq
  defining $t(s)\doteq\int_0^s w_0^q(\sigma)\,d\sigma$, and  $s: [0,+\infty[\to\R_+$ as the (continuous) inverse function of $t(s)$,  the control 
 $\alpha(\cdot)\doteq
 \frac{w(s(\cdot))}{w_0(s(\cdot))}$
 belongs to ${\mathcal A}$ and 
  and $ \xi_{ x}(s(\cdot),w_0,w)$    is the solution of (\ref{S}) corresponding to  $\alpha$.  
 \end{Proposition} 
 
 \begin{Remark}\label{REmb1}{\rm  Considering   extended controls  where $w_0(s)=0$ for  $s$ in some intervals, is a way to introduce a notion of {\it generalized control,} where the (discontinuous) generalized solution   to (\ref{S}) corresponding to $(w_0,w)$, say  $y_x^{gen}$ is defined as  $y_x^{gen}(\cdot)\doteq\xi_x(s(\cdot),w_0, w)$, where  $s(\cdot)$ is, e.g.,  the right inverse  of $t(s)\doteq\int_0^s w_0^q(\sigma)\,d\sigma$ for $s\ge0$. It is clear that, for $q>p$, one has  $f^{\infty}\equiv0$ and    $y_x^{gen}(\cdot)\equiv y_x(\cdot)$   (for more details, see \cite{RS}).}
 \end{Remark}
 
For any $t\ge0$, $x\in\R^n$, we define the {\it extended finite horizon value function}  
$$
 {{ V}} (t,{x})\doteq 
 \inf_{\{(w_0,w)\in\Gamma:  \ \exists S>0 \ \text{s.t.} \ \int_0^Sw_0^q(s)\,ds=t\}}J(S,x,w_0,w)  
$$
and the  {\it extended infinite horizon value function}  
$$
 {{ V}} ({x})\doteq 
 \inf_{(w_0,w)\in\Gamma}J(+\infty,x,w_0,w)  \quad(\le+\infty). 
$$

\begin{Remark}\label{ECset}{\rm In Proposition \ref{REmb}, we establish a correspondence between   $\alpha\in{\mathcal A}$ and   $(w_0,w)\in\Gamma^+$,   assuming (\ref{w-0}). 
This is not a restriction, however, since (\ref{w-0}) is satisfied  by all $(w_0,w)\in\Gamma$ such that $J(+\infty,x ,w_0,w)<+\infty$,  owing to the  coercivity hypothesis (\ref{CC}) which, in the extended problem,  reads as
\begin{equation}\label{CCE}
\bar l(x,w_0,w)\ge C_2|w|^q-C_1w_0^q \qquad \forall (x,w_0,w)\in\R^n\times S(A).
\end{equation}
 In fact,    if we had 
$\int_0^{+\infty}w_0^q(s)\,ds=T<+\infty$,    (\ref{CCE})   together with the constraint $w_0^q+|w|^q= 1$  would yield a cost
$$
J(+\infty, x, w_0,w)\ge C_2\int_0^{+\infty}|w(s)|^q\,ds-C_1 T=+\infty,
$$
which is a contradiction.  

For this reason in the definition of $V(x)$ we can disregard the constraint (\ref{w-0}), which should be naturally assumed,   as in the definition  of $V(t,x)$. This is a key point:  due  to the coercivity hypothesis,  the extended infinite horizon problem reduces to an {\it unconstrained} problem with a {\it compact} control  set. 
}
\end{Remark} 

 In view of  Proposition \ref{REmb} and Remark \ref{ECset}, in the extended setting we can recover ${\mathcal  V}(x)$ and  ${\mathcal  V}(t,x)$  by restricting the minimization to $\Gamma^+$  in the definition of $V(x)$ and $V(t,x)$, respectively.    
  In general,  ${\mathcal  V}(x)$ is neither l.s.c. nor u.s.c.. Moreover, as shown in the following example, if $q=p$    it may happen that  $V(x)<{\mathcal  V}(x)$ at some $x$.
 
 \begin{Example}\label{1}{\rm   
 Let us  consider   the  bi-dimensional control system 
$$
\left\{\begin{array}{l}\dot y_1(t)=\alpha(t) \\
\dot y_2(t)=|y_1(t)|+|y_2(t)|
\end{array}\right.
$$
with  $ y(0)=(y_1(0),y_2(0))=x\in\R^2$  and  $\alpha\in L^1_{loc}(\R_+,\R),$  
and define the cost function
$$
 J(t,x,\alpha)=\int_0^{t}(|y(\tau)|^2+|\alpha(\tau)|)\,d\tau.
 $$
 Since any trajectory issuing from $(1,0)$  has a second component strictly increasing, we get ${\mathcal V}(1,0)=+\infty.$

\noindent  Let us now consider the associated extended system, given by
$$
\left\{\begin{array}{l}
\dot \xi_1(s)=w(s)\\
\dot \xi_2(s)=(|\xi_1(s)|+|\xi_2(s)|)w_0(s),\end{array}\right.
$$
$\xi(0)=(\xi_1(0),\xi_2(0))=x$,  and the extended cost  
$$\begin{array}{l}{J}(S,x,w_0,w)=\int_0^S(|\xi(s)|^2w_0(s)+|w(s)|)\,ds.
\end{array}
$$
Implementing the control $w\doteq-1\chi_{[0,1]}$ the  trajectory  issuing from
$(1,0)$, in time $S=1$ reaches the origin, which is an equilibrium point for the extended system, and the corresponding extended cost is 
$$\begin{array}{l}{J}(+\infty,(1,0),w_0,w)=\int_0^{+\infty}|\xi(s)|^2w_0(s)+|w(s)|\,ds=\int_0^1|w(s)|\,ds=1.
\end{array}
$$
This yields $V(1,0)\le 1$, obviously smaller than ${\mathcal V}(1,0)=+\infty.$
} \end{Example}

When $q=p$, we can prove that   ${\mathcal  V}(x)\equiv V(x)$ using (H2) and the following condition.
  {\it \begin{itemize}
\item[{\bf (H3)}]  Assume that there is some closed set  $\C\subset\R^n$ with compact boundary  such that
 for any $x$ with  $V(x)<+\infty$, there is some $\varepsilon>0$   for which 
\bel{H3}
 \liminf_{s\to+\infty}{\bf d}(\xi_x(s,w_0,w))=0 \quad\text{for any $\varepsilon$-optimal control $(w_0,w)\in\Gamma$.}
 \footnote{Both (H3) and (SC1) below  will also be used in the sequel  for other results  and for $A$ compact. In such a case some obvious changes have to be made ($(w_0,w)\in S(A)$, $V$, and $\bar f$ have to be replaced by $a\in A$, ${\mathcal V}$ and $f$ respectively).}
 \eeq 
\end{itemize}} 

\noindent When ${\mathcal V}\equiv 0$ in $\C$, both (SC1) and  (SC2) below  imply  (\ref{H3}).
{\it \begin{itemize}
\item[{\bf (SC1)}]   There exists a function $U: \R^n\setminus\overset{\circ} \C  \to\R_+$, \,  $C^1$ in $\R^n\setminus\overset{\circ}\C$,  positive definite,  proper on $\C^c$,  such that   $\forall x\in\C^c$,
\begin{equation}\label{C1s}
\max_{ (w_0,w)\in S(A)} \left\{   \langle \nabla U(x), \overline{f}(x,w_0,w)\rangle \right\}\le -  m({\bf d}(x)) 
\end{equation}
for some continuous, increasing function  $m:]0,+\infty[\to]0,+\infty[$.
\end{itemize}}
 
{\it \begin{itemize}
\item[{\bf (SC2)}] There is   some continuous, increasing  function  $c_1:]0,+\infty[\to]0,+\infty[$  such that
\bel{H3suf1}
l(x,a)\ge c_1({\bf d}(x))   \qquad \forall (x,a)\in\C^c\times A.
\eeq 
\end{itemize}}

\noindent  (SC1) means  that  (\ref{SE}) is  UGAS (uniformly globally asymptotically stable) w.r.t. $\partial\C$, so that   {\it all } extended trajectories approach $\C$, at least asymptotically, for any $x\in\C^c$  (see e.g.   \cite{BaRo}).  We point out that (SC1) allows the Lagrangian to be zero outside $\C$.

\noindent  (SC2) instead, involving just the Lagrangian, implies that  $l$ is strictly positive outside $\C$.  For $\C\equiv\{0\}$,  it is  satisfied in LQR problems, where $l(x,a)=x^TQx+a^TRa$ and the matrices $Q$ and $R$ are  symmetric and positive definite.  
 (SC2)  easily implies that $J(+\infty,x, w_0,w)=+\infty$ for any control $(w_0,w)$ not satisfying the $\liminf$-condition in (\ref{H3}),  in view of  Remark \ref{ECset}.  
 
 \vv
We have the following well posedness results.
 \begin{Theorem}\label{extT} For any $t\ge 0$ and $x\in\R^n$, one has

(i) ${\mathcal  V}(t,x)=V(t,x)$ and it is continuous;

(ii) if either $q>p$ or {\rm (H2)} and {\rm (H3)} hold for the same $\C$, then   ${\mathcal  V}(x)=V(x)$.
\end{Theorem}  
 
\noindent {\it Proof.}  Theorem 3.3 in \cite{RS} yields (i) while  Proposition 3.4 in \cite{M} implies (ii) for $q>p$.
It remains to prove  thesis (ii) in case $q=p$. Being $V\le{\mathcal V}$, for any $x\in\C$ the equality  ${\mathcal V}(x)=V(x)=0$ follows trivially  from (H2). Let $x\in\C^c$ and $V(x)<+\infty$ (if $V(x)=+\infty$,  ${\mathcal V}(x)=+\infty$ too). Assume by contradiction that there is some $\eta>0$ such that
$$
V(x)< {\mathcal V}(x)-3\eta.
$$
By hypothesis (\ref{contF}), ${\mathcal V}$ is continuous on  the compact set $\partial\C$, therefore 
\bel{hpC}
{\mathcal V}(\bar x)\le\eta \qquad \forall \bar x\in \C^c  \text{\ such that \ }  {\bf d}(\bar x)<3\delta,
\eeq
for some $\delta>0$. Owing to (H3), there is some   $(\tilde w_0,\tilde w)\in\Gamma$ such that 
$$
\int_0^{+\infty}\ol{l}(\xi_x(s,\tilde w_0,\tilde w),\tilde w_0(s),\tilde w(s))\,ds\le V(x)+\eta
$$
and 
$$
 \liminf_{s\to+\infty}{\bf d}(\xi_x(s,\tilde w_0,\tilde w))=0.
 $$
 Hence, for some $S>0$, we have ${\bf d}(\xi_x(S,\tilde w_0,\tilde w))<\delta$ and,   using the Gronwall's Lemma,  by standard calculations we get that the control $(\tilde w_0^n,\tilde w^n)\in\Gamma^+$ where  $\tilde w^n\doteq \frac{n}{n+1}\tilde w$,   for $n$ large enough  satisfies both  ${\bf d}(\xi_x(S,\tilde w^n_0,\tilde w^n))<2\delta$ and
$$
\int_0^S\ol{l}(\xi_x(s,\tilde w^n_0,\tilde w^n),\tilde w^n_0(s),\tilde w^n(s))\,ds\le \int_0^{S}\ol{l}(\xi_x(s,\tilde w_0,\tilde w),\tilde w_0(s),\tilde w(s))\,ds +\eta.
$$
Thanks to Proposition \ref{REmb}, setting $T\doteq  \int_0^S (\tilde w_0^n)^q(s)\,ds$,  $\exists\, \tilde\alpha\in{\mathcal A}$
corresponding to $(\tilde w_0^n,\tilde w^n)$   such that ${\bf d}(y_x(T,\tilde \alpha))<2\delta$ and 
$$
\int_0^Tl(y_x(t,\tilde\alpha),\tilde\alpha(t))\,dt=\int_0^S\ol{l}(\xi_x(s,\tilde w^n_0,\tilde w^n),\tilde w^n_0(s),\tilde w^n(s))\,ds.
$$
By (\ref{hpC}) it follows that, if $\tilde x\doteq y_x(T,\tilde \alpha)$,  there exists a control $\hat\alpha\in{\mathcal A}$ such that 
$$
\int_0^{+\infty}l(y_{\tilde x}(t,\hat\alpha),\hat\alpha(t))\,dt<\eta. 
$$
Thus the control $\alpha(t)\doteq \tilde\alpha(t)\chi_{[0,T[}(t)+\hat\alpha(t-T)\chi_{[T,+\infty[}(t)$ belongs to ${\mathcal A}$ and satisfies
$$
\int_0^{+\infty}l(y_{  x}(t, \alpha), \alpha(t))\,dt<V(x)+3\eta< {\mathcal V}(x).
$$
At this point the first inequality implies that $ {\mathcal V}(x)<+\infty$, which together with the last inequality yields the required contradiction. Statement (ii) for $q=p$  is therefore proved. \qed
 
 \vv
  ${\mathcal V}(x)$ is in general neither u.s.c. nor l.s.c., even if $A$ is compact.  Sufficient conditions for the upper semicontinuity are given in the following proposition. 
\begin{Proposition}\label{Pu}  Assume that {\rm (H2)} and  {\rm (H3)} hold for the same $\C$.   Then  $Dom({\mathcal V})$ is an open set and ${\mathcal V}$ is locally bounded and u.s.c.  in it.
\end{Proposition}

\noindent{\it Proof .} If $A$ is unbounded  condition (\ref{H3}) is assumed  on the extended trajectories. However, (H2) implies that also in this case (and even if  $q=p$),  for any $x$ with  ${\mathcal V}(x)<+\infty$, there is some $\varepsilon>0$   such that
\bel{H3o}
 \liminf_{t\to+\infty}{\bf d}(y_x(t,\alpha))=0 \quad\text{for any $\varepsilon$-optimal control $\alpha\in{\mathcal A}$.}
\eeq
Indeed, if (\ref{H3o}) were not satisfied for some $x$ and $\alpha$, Proposition \ref{REmb} and the equality ${\mathcal V}(x)=V(x)$ proved in Theorem \ref{extT}, would imply a contradiction:  (\ref{H3}) would not hold for the extended control $(w_0,w)$ corresponding to such an $\alpha$.  From now on,  the proof is the same for  a compact or non compact  set $A$.   

\noindent Fix $\eta>0$ and let $\delta>0$ be as in (\ref{hpC}). Let $x_0\in Dom({\mathcal{ V}})\setminus{\mathcal{ T}}$ and let $\alpha\in{\cal A}$ satisfy 
\begin{equation}\label {c11}
\int_0^{+\infty}l(y_{x_0}(t),\alpha(t))\,dt\le {\mathcal{ V}}(x_0)+\eta,
\end{equation}
where $y_{x_0}(\cdot)\doteq y_{x_0}(\cdot, \alpha)$. In view of (\ref{H3o}) $\exists\,  \bar T$  such that 
$d(y_{x_0}(\bar T))\le\delta\}$. For any $x\in\R^n$, let $y_x(\cdot)\doteq y_x(\cdot,\alpha)$. Estimates (\ref{est1}), (\ref{est2})  imply that  one can choose $\delta'>0$ small enough to have, for all $x\in B(x_0,\delta')$,
\begin{equation}\label {c2}
 |y_x(t)|, \  |y_{x_0}(t)|\le   \bar C , \quad |y_x(t)-y_{x_0}(t)|<\delta''   \qquad  \forall t\in[0,\bar T]
\end{equation}
 for  some $\bar C>0$ and  for any $\delta''>0$.  Now by the Dynamic Programming Principle, in short DPP,  choosing $\delta''\le\delta$, we get
\begin{equation}\label {c3}
\begin{array}{l}
{\mathcal{ V}}(x)\le  \int_0^{\bar T} l(y_x(t),\alpha(t))\,dt+ {\mathcal{ V}}(y_x(\bar T))\le
\int_0^{\bar T} M_{\bar C}(1+|\alpha(t)|^q)\,dt 
+\eta  \le C'
\end{array}
\end{equation}
for some $C'>0$, where the second  inequality holds since ${\bf d}(y_x(\bar T))<2\delta$.
Therefore    $Dom({\mathcal{ V}})$ is  an open set and a simple compactness argument yields that 
$\mathcal{V}$ is bounded on any compact subset of $Dom(\mathcal{ V})$.  
 
The fact that  $\mathcal{ V}$ is u.s.c. in $x_0$  can now be  easily deduced.   Adding and subtracting $\int_0^{\bar T}l(y_{x_0}(t),\alpha(t))\,dt$ to the r.h.s. of (\ref{c3}),  $ \forall x\in B(x_0,\delta')$ one obtains
 $$
\begin{array}{l}
\mathcal{V}(x)\le\int_0^{\bar T} L_{\bar C}(1+|\alpha(t)|^q) |y_x(t)-y_{x_0}(t)|)\,dt+\int_0^{\bar T}l(y_{x_0}(t),\alpha(t))\,dt+\eta\\
 \qquad\qquad\qquad \qquad\qquad\qquad \qquad\qquad\qquad \le L_{\bar C}(\bar T+K)\delta''+\mathcal{V}(x_0)+2\eta,
 \end{array}
$$
where $K\doteq \int_0^{\bar T} |\alpha(t)|^q\,dt$. 
Taking $\delta'$   small enough so that $ L_{\bar C}(\bar T+K)\delta''\le \eta$ one has  
$\mathcal{V}(x)\le\mathcal{V}(x_0)+3\eta,$
and with this the upper semicontinuity of $\mathcal{ V}$ is proved. \qed

 \vv
 Let us observe that the continuity on $\partial\C$ prescribed in (H2) plus (H3) does not yield the lower semicontinuity of ${\mathcal V}(x)$. The continuity of  ${\mathcal V}$ in its whole domain will be discussed in Remark \ref{BCs}.

\subsection{Relaxed problems} \label{Rel}
In this section we introduce the relaxed finite and  infinite horizon problems, for the original problems  when $A$ is compact, and  for the extended problems otherwise.  In order to simplify the notation, the corresponding relaxed value functions,   ${\mathcal V}^r$ (if $A$ is compact) and $V^r$ (in which $A$ is replaced by $S(A)$ and the extended data are considered), will be always denoted by $V^r$.
\vv 
 \noindent {\sc   $A$ compact}.  As  usual we define  the  relaxed   controls  
 $$
 \mu(\cdot)\in {\mathcal A}^r\doteq L^{\infty}(\R_+,{\cal P}(A)),
 $$
 where  $A^r\doteq{\cal P}(A)$  is the set of Radon probability measures on the compact set $A$  endowed with the weak$^*$-topology,  and we consider  
 $\psi\in\{f, \,l\}$  extended to $\R^n\times A^r$ by setting  
$$
 \psi^r(x,\mu)\doteq \int_{A}\psi(x, a)\,d\mu \qquad \forall \mu\in A^r.
 $$
For any $x\in \R^n$ and $\mu\in {\mathcal A}^r$, $y^r_{x}(\t,\mu)$ denotes the relaxed trajectory, solution of 
\bel{SEr}
\dot y^r=f^r(y^r,\mu) \qquad \text{for $\t>0$,} \quad  y^r(0)=x.
\eeq
Finally, we introduce
$$
 V^{r}(t,x)\doteq \inf_{\mu\in {\mathcal A}^r}{\mathcal J}^r(t,x, \mu) \qquad \forall (t,x)\in]0,+\infty[\times \R^n
 $$
 and 
 $$
V^{r}(x)\doteq \inf_{\mu\in {\mathcal A}^r}{\mathcal J}^r(+\infty,x, \mu) \qquad \forall x\in\R^n,
 $$
where
$$
{\mathcal J}^r(t,x, \mu)\doteq\int_0^{t} l^r(y_x^r(\t,\mu),\mu(\t))\,d\t \quad\text{for any $t\in]0,+\infty]$.}
$$
Since  for $A$ compact, 
\bel{con}
\forall x\in\R^n: \qquad \overline{co}(f(x,A)\times l(x,A))= f^r(x,A^r)\times l^r(x,A^r),
\eeq
  standard arguments yield  that the relaxed finite and infinite horizon problems coincide with the original ones under the following convexity hypothesis.

{\it \begin{itemize}
\item[{\bf (CV)}] Let $A$ be compact. For each  $x\in\R^n$,   the  following set is convex:
\bel{convhp}
{\mathcal L}(x)\doteq\left\{ (\lambda,\gamma)\in\R^{n+1}: \      \exists a\in A \ \text{s. t.}  
\ \  \lambda=f(x,a), \  l(x,a)\le\gamma\right\}. 
\eeq
\end{itemize}}
 
\vsm 
\noindent {\sc  $A$  unbounded}. We define relaxed extended controls, 
$$
\mu(\cdot)\in \Gamma^r\doteq L^{\infty}(\R_+,{\cal P}(\overline{B(0,1)}\cap A),
$$
  $A^r\doteq{\cal P}(\overline{B(0,1)}\cap A)$  denotes now the set of Radon probability measures on the compact set $\overline{B(0,1)}\cap A$  endowed with the weak$^*$-topology and we consider  
 $\psi\in\{\overline{f}, \, \overline{l}\}$  extended to $\R^n\times A^r$ by setting  
$$
 \psi^r(x,\mu)\doteq \int_{\overline{B(0,1)}\cap A}\psi(x,(1-|w|^q)^{1/q},w)\,d\mu \qquad \forall \mu\in A^r.
 $$
For any $x\in\C^c$ and $\mu\in\Gamma^r$, $\xi^r_{x}(s,\mu)$ is the relaxed trajectory, solution of 
\bel{SEr}
\dot\xi^r=\overline{f}^r(\xi^r,\mu) \qquad \text{for $s>0$,} \quad  \xi^r(0)=x.
\eeq
In this case, $V^{r}(t,x)$ and $V^{r}(x)$ are given respectively by
$$
 V^{r}(t,x)\doteq \inf_{\{\mu\in\Gamma^r, \, \int_0^S(1-|\mu(s)|^q)\,ds=t\}}{J}^r(S,x, \mu)
 $$
and
 $$
 V^{r}(x)\doteq \inf_{\mu\in\Gamma^r}{J}^r(+\infty,x, \mu),
 $$
 where
$$
{J}^r(S,x, \mu)\doteq\int_0^{S}\overline{l}^r(\xi_x^r(s,\mu),\mu(s))\,ds \quad\text{for any $S\in]0,+\infty]$.}
$$
  
 If $A$ is unbounded,  in order to have   $V^r\equiv V$ we could again invoke  a convexity condition analogous to (CV),   for the extended problem. However, in view of the definitions of $\ol{f}$ and $\ol{l}$  this condition  would be very difficult to be satisfied,  since the control set $S(A)$ is not convex.  Hence we introduce  the  weaker convexity condition (CV)$'$ below, where   $S(A)$ is replaced by $[0,1]\times \left(\overline{B(0,1)}\cap A\right)$ and the space-time extended dynamics $(w_0^q,\ol{f})$ is considered.    (CV)$'$ is  verified,  for instance,  by a control-affine dynamics and  a convex Lagrangian.  

 {\it \begin{itemize}
\item[{\bf (CV)$'$}] Let $A$ be a unbounded.    For any $x\in\C^c$,    the following  set is convex:   
\bel{convhpe1}
\begin{array}{l}
{L}(x)\doteq \left\{(\lambda_0,\lambda,\gamma)\in\R^{1+n+1}: \      \exists (w_0,w)\in[0,1]\times   \left(\overline{B(0,1)}\cap A\right) , \right. \\  \left. \qquad\qquad \text{s.t.}  
 \ (\lambda_0,\lambda)=(w_0^q,\overline{f}(x, w_0,w)), \ \   \overline{l}(x,w_0,w)\le\gamma   \right\}.
\end{array}
\eeq
\end{itemize}} 
\vv 
Both for bounded and unbounded controls,   the relaxed and the original finite horizon problems coincide.   
  \begin{Theorem}\label{Vt,r} {\sc Finite horizon.} For any $(t,x)\in]0,+\infty[\times\R^n$ we have that $V^r(t,x)$ is continuous,  there exists an optimal relaxed control, and
  $$
  \mathcal{V}(t,x)\equiv V^r(t,x).
  $$ 
Moreover,  assuming either {\rm (CV)} or  {\rm (CV)$'$}, there exists an optimal control $\alpha$ for the original problem in case either $A$ is compact or $q>p$, and there exists an optimal extended control $(w_0,w)$  for $p=q$. 
  \end{Theorem}
  
\noindent{\it Proof.} The equality, which could be proved directly,  is a straightforward consequence of   the uniqueness result in Theorem \ref{P1b}, since it is easy to show that  $ V^r(t,x)$ satisfies (\ref{CPb}) in the viscosity sense.  Moreover, it is continuous as ${\mathcal V}(t,x)$, since the relaxed data have the same properties of the original ones.   The  existence of an optimal control for the relaxed problem (which does not imply in general the existence of an optimal ordinary  control) is well known.

If  (CV)   holds, an optimal control $\alpha$ for ${\mathcal V}(t,x)$ exists by standard arguments. When $A$ is unbounded,   in view of (CV)$'$, in correspondence to an optimal relaxed control $\mu_r$ for $V^r(t,x)$,   there is a control $(w_0,w)\in{\mathcal B}(\R_+, [0,1]\times   (\overline{B(0,1)}\cap A ))$ such that 
$\xi(\cdot)\doteq \xi^r_x(\cdot,\mu)\equiv \xi_x(\cdot,w_0,w)$,    $J^r(S, x, \mu)\ge J(S, x, w_0,w)$ and in addition
\bel{002}
\int_0^{S}w_0^q(s)\,ds=\int_0^{S}(1-|\mu(s)|^q)\,ds=t
\eeq
for some $S>0$. In general,  $(w_0,w)\notin \Gamma$ since $w_0^q+|w|^q$ may differ from 1. Nevertheless,  using the arc-lenght reparameterization $\Phi^{-1}$, where 
$\Phi(\sigma)=\int_0^{\sigma}[w_0^q(s)+|w(s)|^q]\,ds$,  
 the control $(w_0,w)$ can be substituted by one taking values in $S(A)$,  satisfying (\ref{002}), and having the same cost and trajectory. This is possible since  $\overline{f}$ and $\overline{l}$ are $q$-positively  homogeneous in $(w_0,w)$ (see also Proposition \ref{REmb}).  Such a control is clearly the desired optimal extended control.     

\noindent When $q>p$,  we show that, in correspondence to any extended control $(w_0,w)\in\Gamma$ verifying (\ref{002}) and $J(S, x,w_0,w)<+\infty$,  there exists  $\alpha\in{\mathcal A}$ such that 
$$
{\mathcal J}(t, x,\alpha)\le J(S, x,w_0,w).
$$
  Suppose first that   $w_0=0$ on a unique (bounded) interval $[s_1, s_2]$. Then the   trajectory $\xi_x(s, w_0,w)\equiv \xi_x(s_1, w_0,w)$  for all $s\in [s_1, s_2]$   because of the definition of $f^\infty$, while $l\ge0$ implies that  $\int_{s_1}^{s_2}\overline{l}(\xi_x(s, w_0,w), w_0,w)\,ds\ge 0$.   Therefore   $J(S,x,w_0,w)\ge J(S-(s_2-s_1),x, \tilde w_0,\tilde w)$ if $(\tilde w_0,\tilde w)(s)\doteq \chi_{[0,s_1[}(w_0,w)(s)+\chi_{[s_1,S-(s_2-s_1)]}(w_0,w)(s+s_2-s_1)$ for all $s\in[0,S-(s_2-s_1)]$.  For the general case, set $\sigma=\sigma(s)\doteq \int_0^s \chi_{]0,1]}(w_0(s'))\,ds'$ and let $s=s(\sigma)$ be the right inverse of $\sigma(\cdot)$.  It is easy to see that the  control  $(\tilde w_0,\tilde w)(\sigma)\doteq (w_0,w)(s(\sigma))$ for all $\sigma\ge0$ does the job. 
 The above argument lets us immediately  conclude in view of Proposition \ref{REmb}, since $(\tilde w_0,\tilde w)\in\Gamma^+$.
\qed
  \vv
As it is well known, this relaxation property is no more true for the infinite horizon problem and ${\mathcal V}(x)$  does not coincide in general with  $V^r(x)$,  even in the simplest case of compact valued  controls, as shown by Example \ref{ExC} below. The following weaker results hold.
 \begin{Theorem}\label{COVr} {\sc Infinite Horizon.} 
 \begin{itemize}
\item[(i)] Assume either  {\rm (CV)} or {\rm (CV)$'$} and $q>p$.  Then for any $x\in\R^n$ we have  
\bel{=}
{\mathcal V}(x)=V^{r}(x)
\eeq
and  there exists an optimal control $\alpha\in{\mathcal A}$ for the original   problem.
\item[(ii)]  Assume   {\rm (CV)$'$} and $q=p$. Then for any $x\in\R^n$,  
\bel{=e}
V(x)= V^r(x)
\eeq
  and there exists an optimal extended control, $(w_0,w)\in\Gamma$. If  moreover {\rm (H2)} and  {\rm (H3)}  hold for the same $\C$,  then we have {\rm (\ref{=})}.
\end{itemize}
\end{Theorem}
\noindent {\it Proof.} Let us   prove that, assuming (CV)$'$,  $V(x)=V^r(x)$ for $q\ge p$.
Let  $x\in\R^n$ be such that $V^r(x)<+\infty$ (if $V^r(x)=+\infty$, $V(x)=+\infty$ too).   In order to prove  (\ref{=e}), let $\mu\in\Gamma^r$ be an optimal  relaxed control,
such that 
$$
J^r(+\infty, x, \mu)\doteq \int_0^{+\infty}\overline{l}^r(\xi_x^r(s,\mu),\mu)\,ds=V^r(x)<+\infty, 
$$ 
  whose existence is proved in Theorem \ref{supC} below. Thanks to (CV)$'$, by standard arguments there exists a control $(w_0,w)\in{\mathcal B}(\R_+, [0,1]\times   (\overline{B(0,1)}\cap A))$ such that 
$\xi^r_x(\cdot,\mu)\equiv \xi_x(\cdot,w_0,w)$,    $J^r(+\infty, x, \mu)\ge J(+\infty, x, w_0,w)$ and
$$
\int_0^{\sigma}w_0^q(s)\,ds=\int_0^{\sigma}(1-|\mu(s)|^q)\,ds \qquad \forall \sigma\ge0.
$$
From the same arguments  in Remark \ref{ECset} applied to the relaxed problem, we have that
$
\int_0^{+\infty}w_0^q(s)\,ds=+\infty.
$
Now, $(w_0,w)\notin \Gamma$  in general, but by using the arc-lenght reparametrization and arguing as in the proof of  Theorem \ref{Vt,r},  we can obtain an  extended control in $ \Gamma$ with the same cost,  and this proves  (\ref{=e}).  The last statement of (ii) follows from Theorem \ref{extT} (ii).
 
\noindent If $A$ is compact, statement (i) can be proved by standard arguments. When $A$ is unbounded, the equality ${\mathcal V}(x)=V^{r}(x)$  follows from the previous point together with Theorem \ref{extT} (ii).  The existence of an optimal control $\alpha$ in the case $q>p$ can be recovered as in the last part of the proof of Theorem \ref{Vt,r}.  \qed

\begin{Remark}{\rm In case $A$ unbounded and $q=p$, even if ${\mathcal V}\equiv V^r$, both  the original finite and  infinite horizon problems may not have an optimal control.}
\end{Remark}
 
 \section{Finite-horizon  approximation } \label{FHA}
In this section 
we give a representation formula for the limit, as $t$ tends to $+\infty$ of the finite horizon value functions 
$$
{\mathcal V}(t, x)\doteq\inf_{\alpha\in{\mathcal A} }\int_0^tl(y(\t),\alpha(\t))\,d\t,
$$
defined as
\bel{sig}
\Sigma(x)\doteq\lim_{t\to+\infty}{\mathcal V} (t,x)=\sup_{t>0}{\mathcal V} (t,x) \qquad\forall x\in\R^n.
\eeq

The following simple example describes what is expected to happen, for the compact control case.  
 \begin{Example}\label{ExC}{\rm 
Let us  consider   the  bi-dimensional control system 
$$
\left\{\begin{array}{l}
\dot y_1=\alpha(t) \\
\dot y_2=|y_1(t)|
\end{array}\right.
$$
with  $ y(0)=x\in\R^2,$  $\alpha(t)\in A\doteq\{\pm1\},$ and define
 the cost function
$$
 J(t,x,\alpha)=\int_0^{t}|y(\tau)|^2\,d\tau.
 $$
Clearly, any trajectory issuing from $(0,0)$ has a strictly increasing second component, which gives immediately
 ${\mathcal V}(0,0)=+\infty,$ while the relaxed value function $V^r(0,0)=0$. $V^r$, indeed, coincides with the infinite horizon value function where controls $\alpha(t)\in[-1,1]$ are allowed. 
 
  \noindent 
Now fix $t>0$ and for every $n\in \N$, $n>0$ let us set $h\doteq\frac{t}{n}$ and let us define the control
 $$\alpha_n\doteq(-1)^i\quad\forall \tau\in[ih,(i+1)h),\quad i=0,\dots,n-1.$$ The trajectory issuing from $(0,0)$, relative to $\alpha_n,$  
 has the first component  such that $\sup_{[0,t]}|y_1(t,\alpha_n)|\le\frac{t}{n}$ and for the second component
 $\sup_{[0,t]}|y_2(t,\alpha_n)|\le\frac{t^2}{n}$ which gives
$$
 J(t,x,\alpha_n)=\int_0^{t}|y(\tau)|^2\,d\tau\le \frac{t^3(1+t^2)}{n^2},
 $$
 and this yields  ${\mathcal V}(t,(0,0))=0$  for every $t>0$. Therefore, $\Sigma(0,0)=0=V^r(0,0)$.
   }
 \end{Example}  
 
 The result suggested by the previous example can be extended to the case of unbounded controls as follows.
\begin{Theorem}\label{supC} 
For any $x\in\R^n$,  we have
$$
\Sigma(x) =V^r(x).
$$
Moreover, $V^r$ is l.s.c. and there exists an optimal relaxed control.
\newline
\end{Theorem}
In case $A$ unbounded, we use the following preliminary result,  true thanks to the coercivity hypothesis (\ref{CC}) and interesting  in itself.
\begin{Proposition}\label{underline V} For any $x\in\R^n$, 
$$
\Sigma(x)= \sup_{s>0}W(s,x),
$$
where 
$$
W(s,x)\doteq \inf_{\mu\in\Gamma^r}\int_0^{s}\overline{l}^r(\xi_x^r(s,\mu),\mu(s))\,ds.
 $$
\end{Proposition}
\noindent {\it Proof. }  Let $x\in\R^n$. 
  We recall that for any $t>0$,   $\mathcal{V}(t,x)$ coincides with the relaxed finite horizon value function $V^r(t,x)$ in view of Theorem \ref{Vt,r}. Hence  $\Sigma(x)=\sup_{t>0}V^r(t,x).$ In order to conclude, it remains essentially  to prove   that the time constraint  $\int_0^S(1-|\mu(s)|^q)\,ds=t$ in the definition of  $V^r(t,x)$ can be dropped, so that  
$$
\sup_{t>0}V^r(t,x)= \sup_{s>0}W(s,x).
$$
 Let us  first show the simpler inequality 
 \bel{in1}
\Sigma(x)\ge \sup_{s>0}W(s,x),
 \eeq
 true even in non coercive problems.  By Theorem \ref{Vt,r}, for any $n\in\N$, there exists an optimal relaxed trajectory-control pair $(\xi^r_n,\mu_n)$ and some $s_n>0$ such that 
 $$
 V^r(n,x)=\int_0^{s_n}\overline{l}^r(\xi_n^r(s,\mu_n),\mu_n(s))\,ds, \qquad \int_0^{s_n}(1-|\mu_n(s)|^q)\,ds=n.
 $$
Hence 
\bel{in2}
V^r(n,x)\ge W(s_n,x)
\eeq
where $s_n\ge n$ by definition,  so that (\ref{in1}) follows easily by passing to the limit as $n$  tends to $+\infty$  in (\ref{in2}) (the $\lim_{s\to+\infty} W(s,x)$  exists and coincides with $\sup_{s>0}W(s,x)$ by monotonicity).  
 
 \noindent Now, by (\ref{in1}) the converse inequality is trivially satisfied if $\sup_{s>0}W(s,x)=+\infty$. Let us assume  by contradiction that there is some $\eta>0$ such that
 \bel{in3}
  \sup_{s>0}W(s,x)<  \Sigma(x)-\eta.
 \eeq
 Then for any $n\in\N$ there is some $(\xi^r_n,\mu_n)$ such that
$$
\int_0^{n}\overline{l}^r(\xi_n^r(s,\mu_n),\mu_n(s))\,ds <  \Sigma(x)-\eta.
$$
 Let $t_n\doteq \int_0^{n}(1-|\mu_n(s)|^q)\,ds$ \,($\ge0$). If $\{t_n\}_n$ is unbounded,   for some subsequence, still denoted by $\{t_n\}_n$, $t_n>0$ for all $n$,  $\lim_n t_n=+\infty$ and we get
 $$
 V^r(t_n,x)\le \int_0^{n}\overline{l}^r(\xi_n^r(s,\mu_n),\mu_n(s))\,ds <  \Sigma(x)-\eta.
  $$
 Thus   letting $n$ tend to $+\infty$  one obtains that $\Sigma(x)=\lim_n V^r(t_n,x)\le \Sigma(x)-\eta$, which yields the  desired contradiction. 
 
\noindent If instead the sequence  $\{t_n\}_n$ is bounded, so that  $t_n\le T$ for all $n$ for some $T>0$ by  the coercivity assumption (\ref{CC}) we get  
$$
C_2n-(C_2+C_1)T\le C_2\int_0^{n}|\mu_n(s)|^q\,ds-C_1\int_0^{n}(1-|\mu_n(s)|^q)\,ds< \sup_{s>0}W(s,x)<+\infty.
$$
When $n$ tends to $+\infty$, the l.h.s. tends to  $+\infty$ and we  get   a contradiction also in this case. \qed

\vv
\noindent {\it Proof of Theorem \ref{supC}. } We consider only the case $A$ unbounded, the proof for $A$ compact being similar and actually simpler.
 By the previous proposition, 
$\Sigma(x)= \sup_{s>0}W(s,x)\le V^r(x),$   being $\ol{l}\ge0$.     
When $\Sigma(x)=+\infty$, we have trivially $\Sigma(x)=V^r(x).$ Let  thus suppose $\Sigma(x)<+\infty$. For every $n\in\N$ there exists an optimal  relaxed trajectory-control pair $(\xi^r_n,\mu_n)$ satisfying
\begin{equation}\label{c1}
\Sigma(x)=\lim_n W(n,x)=\lim_n\int_0^{n}\overline{l}^r(\xi_n^r(s),\mu_n(s))\,ds.
\end{equation}
Let $S>0$.  Owing to the compactness of the control set $\overline{B(0,1)}\cap A$,   the set $\{\xi_n^r\}_n$ is uniformly bounded and equilipschitz on $[0,S]$. Moreover, for any $n\ge S$, 
$$
\int_0^{S}\overline{l}^r(\xi_n^r(s),\mu_n(s))\,ds\le \Sigma(x).
$$ 
Therefore by Ascoli-Arzel\`a Theorem  there exists a subsequence  $\{\xi^r_{n'}\}_{n'}$,   uniformly converging to some function $\bar \xi^r$ in $[0,S]$, such that, owing to (H0),   
\bel{eqlim}
\int_0^{S}\overline{l}^r(\bar\xi^r(s),\mu_n(s))\,ds\le \Sigma(x)+\rho_S(n),
\eeq 
for some $\rho_S(n)$ with $\lim_n\rho_S(n)=0$. 
Moreover, since $L^{\infty}([0,S],{\cal P}(\overline{B(0,1)}\cap A)$ is sequentially weakly$^*$-- compact (see \cite{W},  p. 272), there exists a subsequence  $\{\mu_{n''}\}_{n''}$ of $\{\mu_{n'}\}_{n'}$  which converges weakly to some $\bar\mu$ in $[0,S]$. 
 Therefore  by a diagonal procedure  we  obtain a trajectory-control pair $(\bar \xi^r,\bar\mu)$ defined on the whole interval $\R_+$ and such that for any $S>0$ there is some subsequence $\{(\xi_n^r,\mu_n)\}_n$,  where $\xi_n^r$ converges uniformly to $\bar\xi^r$ and $\mu_n$ weakly to $\bar\mu$ in $[0,S]$. 
 
\noindent For any $S>0$, by the weak convergence, passing to the limit in (\ref{eqlim}) one has
$$
\int_0^{S}\overline{l}^r(\bar\xi^r(s),\bar\mu(s))\,ds\le \Sigma(x).
$$
Consequently, since $\ol{l}$ is nonnegative,  
$
V^r(x)=\int_0^{+\infty}\overline{l}^r(\bar\xi^r(s),\bar\mu(s))\,ds=\Sigma(x)
$
 (and  $\bar\mu$ is   the optimal  relaxed control).

\vv
We are going now to discuss the relation of the previous approximation result with the original value function
${\mathcal V}$. A straightforward consequence of Theorems \ref{COVr} and  \ref{supC} is the following

\begin{Corollary}\label{CorsupC}
Assume  either  {\rm (CV)} or {\rm (CV)$'$}. If $A$ is unbounded and $q=p$ let  {\rm (H2)} and  {\rm (H3)} hold  for the same $\C$.   Then for any $x\in\R^n$ we have  
$$
\Sigma(x)={\mathcal V}(x), 
$$
where $\Sigma$ is defined in {\rm (\ref{sig})}.
\end{Corollary}

If no convexity is assumed, we prove that $\Sigma(x)={\mathcal V}_*(x)$, the l.s.c. envelope of  ${\mathcal V}$, under  some mild additional hypotheses  (H0)$_1$  and  (H0)$_2$.   Let us remark that, since the boundary value problem associated to the infinite horizon value function considered here has not a unique solution, we have to prove this relaxation result directly. 

{\it \begin{itemize}
\item[${\bf (H0)_1}$]  (i) Hypothesis {\rm (H0)} holds with  the constants $L_R$, $M_R>0$ and the modulus  $\omega(\cdot)\doteq\omega(\cdot,R)$ independent of $R$ and  
$$
|f(x,a)|\le M(1+|a|^p)\qquad\forall x\in\R^n, \ a\in A.
$$
(ii) Moreover, $\int_0^1(\omega(s)/s)\,ds<+\infty$.
\end{itemize}}   

{\it \begin{itemize}
\item[${\bf (H0)_2}$]  (i) For every $x\in\R^n$ with $V^r(x)<+\infty$ there exists an {\rm optimal} relaxed  control $\mu$ such that,  for some $\bar R>0$,
\bel{brt}
|\xi_x^r(s,\mu)|\le\bar R \quad \forall s\in[0,+\infty[, 
\eeq
if $A$ is unbounded   $[|y_x^r(t,\mu)|\le\bar R$ \, $\forall t\in[0,+\infty[$,  if $A$ is compact$]$.

 \noindent (ii)  Moreover $\int_0^1(\omega(s, \bar R+3)/s)\,ds<+\infty$,  where $\omega$ is the modulus of $l$ introduced in {\rm (H0)}.
\end{itemize}}

Hypothesis   (H0)$_2$ (i) roughly says that relaxed  trajectories going to infinity are  not convenient. Both hypotheses (SC1) and (SC2) introduced in Section \ref{Gen} yield {\rm { (H0)$_2$}} (i). Actually, we recall that condition (SC1)  implies the UGAS property w.r.t. $\partial\C$ for the {\it relaxed} control system too. Therefore, {\it all}  the relaxed  trajectories approach  the compact  set $\partial\C$ asymptotically (see e.g. \cite{BaRo}). This easily implies {\rm { (H0)$_2$}} (i).  (SC2) instead,  implies (\ref{linfty}) below, which we will show to be  sufficient for  {\rm { (H0)$_2$}} (i) in Proposition \ref{suff(H2)}.  Conditions  (H0)$_1$ (ii) and   (H0)$_2$ (ii) are    fulfilled, e.g.,  if   $\omega(r)=Lr^\gamma$ and $\gamma>0$. 

 \begin{Theorem}\label{Thril}  Assume either  {\rm (H0)$_1$} or  {\rm (H0)$_2$}.   
 \begin{itemize}
 \item[(i)]  If either $A$ is compact or  $q>p$, then for any $x\in\R^n$,
\bel{l.s.c.1}
 {\mathcal V}_*(x)=V^r(x); 
\eeq
  \item[(ii)] if $A$ is unbounded  and $q=p$, then for any $x\in\R^n$,
\bel{l.s.c.2}
V_*(x)=V^r(x). 
\eeq 
Moreover, if {\rm (H2)} and  {\rm (H3)} hold for the same $\C$,  we have  {\rm (\ref{l.s.c.1})}.
 \end{itemize}
  \end{Theorem}
  
 \noindent{\it Proof.}  We prove the theorem only for $A$ unbounded, the proof for $A$ compact being analogous and actually simpler.  We   show    that  (\ref{l.s.c.2}) holds for any $x\in\R^n$.   Both   statement (i) for $q>p$ and the last part of  (ii) for $q=p$  follow then from Theorem \ref{extT} (ii). 
 
 \noindent Since  $V^r(x)\le V(x)$ and $V^r$  is l.s.c.,  then $V^r(x)\le V_*(x)$   for any $x\in\R^n$. It  remains to prove the converse inequality, where it is not restrictive to consider only $x\in\R^n$ with $V^r(x)<+\infty$. 
 
 Let us first assume  (H0)$_1$. In this case it is easy to prove that $\ol{f}$ and $\ol{l}$ verify Assumption 3.1 of \cite{AB}, so that (\ref{l.s.c.2}) holds in view of Theorem 3.2 of the same paper.   Actually, in \cite{AB} infinite horizon problems in $L^\infty$ are considered, but for a nonnegative running cost  $\ol{l}$, one has 
 $$
 \underset{s\in[0,+\infty[}{\text{ess}\sup}\int_0^s\ol{l} (\xi(s),w_0(s),w(s))\,ds=\int_0^{+\infty} \ol{l} (\xi(s),w_0(s),w(s))\,ds.
 $$  
  
 Let now  (H0)$_2$ be in force.  Accordingly, let   $(\xi^r(\cdot),\mu(\cdot))$, where $\xi^r(\cdot)\doteq \xi_x(\cdot,\mu)$,  be a relaxed  optimal trajectory-control pair satisfying
(\ref{brt})  for some $\bar R>0$.  Let $\psi:\R^n\to [0,1]$ be a  $C^\infty$ cut-off map such that for all $x\in\R^n$,
 $$
 \psi(x)=1 \quad \text{if } \ |x|\le \bar R+1; \qquad \psi(x)=0 \quad \text{if } \ |x|\ge \bar R+3.
 $$
  Now  $\ol{f}_{\bar R}\doteq \psi\,\ol{f}$, $\ol{l}_{\bar R}\doteq \psi\,\ol{l}$ satisfy hypothesis  (H0)$_1$ and thus Assumption 3.1 of \cite{AB}.  Hence by the proof of Theorem 3.2 in \cite{AB}, for any $\varepsilon>0$  there exist an extended control $(w_0,w)\in\Gamma$ and an extended trajectory $\xi(\cdot)$  such that $\dot\xi(s)=\ol{f}_{\bar R}(\xi(s),w_0(s),w(s))$ for a.e. $s\in]0,+\infty[$  and
 \bel{AB1}
  |\xi(s)-\xi^r(s)|\le 2\varepsilon\,e^{-2\bar L s} \qquad \text{for} \ s\in]0,+\infty[,
 \eeq
\bel{AB2}
  \int_0^{+\infty}\ol{l}_{\bar R} (\xi(s),w_0(s),w(s))\,ds\le J^r(+\infty, x, \mu)+\varepsilon +\int_0^{2\varepsilon}\frac{\omega(s,\bar R+3)}{2\bar L s}\,ds,
\eeq
  where $\bar L>0$ is the Lipschitz constant of $\ol{f}_{\bar R}$ (which can be assumed equal to $L_{\bar R+3}$) and $\omega$ is the same as in  (H0)$_2$. 
 Set $\bar x\doteq \xi(0)$. From (\ref{AB1}) it follows that  $|\xi(s)|<\bar R+1$ for all $s\ge0$ as soon as $\varepsilon<1/2$. Hence in view of the definition of $\ol{f}_{\bar R}$ and $\ol{l}_{\bar R}$,  $\xi(\cdot)$ solves the original system (\ref{SE}) with initial condition $\bar x$ and (\ref{AB2}) holds with $\ol{l}_{\bar R}$ replaced by $\ol{l}$. Taking the limit as $\varepsilon$ tends to zero we conclude that $V_*(x)\le V^r(x)$. \qed
  
\vv
A sufficient condition to have  (H0)$_2$ (i),  is given in the  
next proposition. Let us remark  that  (\ref{linfty}), even in the case  $A$   unbounded,   involves only the original Lagrangian $l$ and not the extended $\overline l$. 

\begin{Proposition}\label{suff(H2)} Let us assume that, for every $x\in\R^n$,  
\bel{linfty}
\liminf_{|x|\to+\infty}\left(\inf_{a\in A}l(x,a)\right)>0.
\eeq
Then {\rm { (H0)$_2$}} {\rm (i)} holds.
\end{Proposition}

\noindent {\it Proof. }
Let  $A$ be  unbounded. Then condition (\ref{linfty}) together with the coercivity assumption (\ref{CC}) easily implies  
\bel{s1}
\ol{l}(x,w_0,w)\ge \bar C \quad \forall x \ \text{ with $|x|\ge \bar M$ and } \ (w_0,w)\in S(A)
\eeq
for some positive constants $\bar M$,  $\bar C$, so that the same holds true for $\ol{l}^r$.  Assume by contradiction that for some $x$ with $V^r(x)<+\infty$,  there exists some optimal relaxed control $\mu$ such that the corresponding trajectory $\xi^r(\cdot)$ satisfies $|\xi^r(s_n)|\ge n$ for  some increasing, positive sequence $s_n$ tending to $+\infty$. Then  $\exists N>0$ such that $|\xi^r(s_n)|>\bar M$ for all $n\ge N$. If   $|\xi^r(s)|>\bar M$ for all $s\ge s_n$ for some $n$, then by (\ref{s1}) we should have  an infinite cost, while $J^r(+\infty,x,\mu)=V^r(x)<+\infty$.
Otherwise, we can suppose that  for any  $n>N$ there exists $s_{n+1}\ge c_n>s_n$ such that $|\xi^r(s)|>\bar M$ for $s\in[s_n,c_n[$ and $|\xi^r(c_n)|=\bar M$.  Then by the estimate
 $$
 c_n-s_n\ge\frac{1}{M}\log\left(1+\frac{n-\bar M}{1+\bar M}\right)
 $$
proved in  Lemma 1, pag. 778 of \cite{B}, where $M$ is the constant in (\ref{hp0}), we get
$$
\int_0^{+\infty}\ol{l}^r(\xi^r_x(s), \mu^r(s))\,ds\ge \sum_{n=N}^{+\infty} \bar C(c_n-s_n)\ge 
\frac{1}{M}\bar C \sum_{n=N}^{+\infty}\log\left(1+\frac{n-\bar M}{1+\bar M}\right)=+\infty,
$$
that is, the same contradiction as above.

\noindent Wo omit the proof in the case $A$ compact, since it is completely similar. \qed
\vv
 In many applications (\ref{linfty}) holds since for some  $r>0$,  $l$ satisfies the following stronger version of (\ref{CC}):
\bel{CCs}
l(x,a)\ge C_2|a|^q+C_1|x|^r \qquad \forall(x,a)\in\R^n\times A
\eeq 
where $C_1$, $C_2>0$ and $q\ge p$ is the same as in (H0).  Condition (\ref{CCs}) holds,   for instance, for  in LQR problems, where $l(x,a)=x^TQx+a^TRa$ and the matrices $Q$ and $R$ are  symmetric and positive definite.

\section{Maximal and minimal solutions and uniqueness}\label{Un}
In this section we  give sufficient conditions in order to characterize  ${\mathcal V}(x)$ as unique solution of the associated HJB equation introduced below. As a byproduct we also obtain the characterization of  the limit function $\Sigma(x)=V^r(x)$.   We start by recalling a uniqueness theorem for the finite horizon problem obtained in \cite{RS} (see also  \cite{MS2}, where more general results, including second order PDEs, are obtained).   We point out that  these results  cannot be derived by classical theorems within  the viscosity theory, in view of the hypothesis $l\ge0$ and of the growth of the data considered here.  Then we derive from the results in \cite{M} and  \cite{MS}  a  uniqueness theorem for the infinite horizon case, generalizing that  obtained for $A$ compact   in \cite{MS1}. 

\vv
Let us define  the  Hamiltonian
\bel{ham1}
{\mathcal H}(x,p)\doteq\sup_{a\in A}\left\{-\langle f(x,a),\, p\rangle-l(x,a)\right\} \qquad \forall (x,p)\in\R^{2n}. 
\eeq
Notice that  in case $A$ unbounded and  $p=q$,  ${\mathcal H}$ can be discontinuous and equal to $+\infty$ at some points.  
When $A$ is unbounded and $q\ge p$, ${\mathcal H}$ can be replaced, as shown in \cite{RS} and \cite{M}, by  the {\it extended} Hamiltonian 
\bel{ham2}
H(x,p)\doteq\max_{(w_0,w)\in S(A)
}\left\{-\langle\overline{f}(x,w_0,w),\, p\rangle-\overline{l}(x,w_0,w)\right\} \qquad \forall (x,p)\in\R^{2n},
\eeq
which  turns out to be  continuous.  Actually,  considering  $H$ is  useful  even if $q>p$, since it allows to consider dynamics verifying $|f(x,a)|\le M(1+|a|^p)(1+|x|)$ instead of the more restrictive hypothesis  $|f(x,a)|\le M(1+|a|^p+|x|)$, assumed  in most of the literature (see e.g. \cite{BDL}, \cite {DL},   and more recently,   \cite{GSor} and the references therein).  An analogous remark holds for $l$. 
 Therefore in the sequel we will use $H$ and, in order to unify the exposition, we  will set $H\doteq{\mathcal H}$ when $A$ is compact.
 
\begin{Example}{\rm In  control-affine problems, or, more precisely, when   $A$ is  unbounded, $q=p=1,$ and $\forall(x,a)\in\R^n\times A$ we have
$$
f(x,a)=f_0(x)+\sum_{i=1}^{m}f_i(x)a_i, \qquad l(x,a)=l_0(x)+\sum_{i=1}^{m}l_i(x)a_i+l_\infty(x)|a|,
$$
 we showed in Section 5 of \cite{MS2},  that  the evolutive PDE is equivalent to the following quasi-variational inequality: 
$$
\max\left\{u_t-\left\langle f_0(x),\, Du(x)\right\rangle-l_0(x),\, \, K(x,Du(x))-l_\infty(x) \right\}=0,
$$
where
$$
K(x,p)\doteq  \max_{w\in A, |w|=1}\left\{-\left\langle \sum_{i=1}^{m}f_i(x)w_i,\, p\right\rangle -\sum_{i=1}^{m}l_i(x)w_i\right\}.
$$
An analogous equivalence holds for the stationary equation.
This is the more usual  formulation of the PDE associated to impulsive control problems.
 }
 \end{Example}
 
For the finite horizon problem we recall what follows. 
 \begin{Theorem}\label{P1b}  {\rm [Corollary 2.1, RS]} We have  ${\mathcal V}(t,x)=V(t,x)$ and it is  continuous  for any $(t,x)\in\R_+\times \R^n$. Moreover, for every $T>0$, it is the unique viscosity solution of the Cauchy problem 
\bel{CPb}
\left\{\begin{array}{l}
u_t+H(x,Du(x))=0 \qquad \forall (t,x)\in]0,T[\times\R^n \\ \, \\
u(0,x)=0 \qquad \forall x\in\R^n
\end{array}\right.
\eeq
 among the functions bounded from below and continuous on $(\{0\}\times\R^n)\cup (\{T\}\times\R^n)$.
 \end{Theorem}
 
 The above uniqueness result, for the case $A$ compact,  can be found in \cite{BCD}.  For $A$ unbounded,  some  comparison theorems due in\cite{BDL} (for the finite horizon problem)   and in \cite{DL} (for the infinite horizon case),  address  just the coercive case $q>p$, as observed above,  require stronger hypotheses on $f$ and $l$,  and imply uniqueness in the class of  the locally Lipschitz functions. We refer to \cite{G} for a uniqueness result among convex functions.

\vv
Leu us now consider the infinite horizon problem with HJB equation
\bel{Ei}
H(x,Du(x))=0.
\eeq
In order to apply the results of \cite{M},  from now on  we assume that 
$$
\text{{\it for any $R>0$, there exists $\bar L_R>0$ such that  $\omega(r,R)=\bar L_R\,r$,} }
$$
  where $\omega$ is the modulus  of continuity of $l$ in (H0).\footnote{The sublinear growth of $l$ assumed in \cite{M} can be removed   as in \cite{GSor}.} 
We recall
\begin{Theorem}\label{5.2}  {\rm [Theorem 4.5, M]}
(i) \, $V\le u$ for any  nonnegative and continuous   supersolution $u:\R^n\to\R\cup\{+\infty\}$   to {\rm (\ref{Ei})} in $\R^n$;

\noindent (ii)\,  $V^{r}(=\Sigma)$   is   l.s.c and it is  the minimal nonnegative   supersolution   to { \rm   (\ref{Ei})} in $\R^n$.\footnote{A function $u:\R^n\to\R\cup\{+\infty\}$ is a viscosity supersolution to {\rm (\ref{Ei})} at $x$  if  either $u_*(x)=+\infty$ or, if $u_*(x)<+\infty$,  it is a supersolution at  $x$. }
\end{Theorem}
Let us  set
$$
\begin{array}{l}{\mathcal S}\doteq 
\left\{(u,\Omega), \quad \Omega\subset\R^n \ \text{open, and }  u:\R^n\to\R_+\cup\{+\infty\}, 
\text{ supersolution  }\right.\\  
  \quad\qquad\left.  
\text{of  { \rm   (\ref{Ei})} in $\R^n$, locally bounded subsolution of  { \rm   (\ref{Ei})}  in 
 $\Omega$, and}\right. \\ \left. 
\qquad\qquad\qquad\qquad  \lim_{x\to\bar x}u(x)=+\infty \quad \forall \bar x\in\partial \Omega.
\right\}\end{array}
$$
The proof of the following theorem follows from Theorem \ref{ET3} below. 
\begin{Theorem}\label{EP3} Assume  {\rm (H2)} and  {\rm (H3)}  for the same $\C$, and alternatively (i) or (ii) below.   
\begin{itemize}
\item[(i)] Assume that either  {\rm (H0)$_1$} or {\rm (H0)$_2$}  holds. Moreover, let ${\mathcal V}$  be continuous in $Dom({\mathcal V})$ and satisfy the boundary condition
\bel{Bcn}
\lim_{x\to\bar x}{\mathcal V}(x)=+\infty \quad \forall \bar x\in\partial Dom({\mathcal V});
\eeq
\item[(ii)]\, assume that either {\rm (CV)} or {\rm (CV)$'$} holds. 
\end{itemize}  
Then 
$\mathcal{V}$  $(\equiv V^r\equiv \Sigma)$   is the unique nonnegative  viscosity solution to  { \rm   (\ref{Ei})} in $Dom({\mathcal V})$,   among the pairs $(u,  \Omega)$ in ${\mathcal S}$, where  $\Omega\supset{\mathcal T}$,  $u\equiv 0$ on $\mathcal{T}$.  Moreover ${\mathcal V}$ is continuous.

\noindent If we drop  {\rm (H0)$_1$},  {\rm (H0)$_2$} in (i), ${\mathcal V}$ (possibly $\ne V^r$) is the unique solution just  among the continuous functions.  
 \end{Theorem}
 
 By the Kruzkov transform $\Psi(v)\doteq 1-e^{-v}$, the above free boundary problem,
   can be replaced by  another boundary value problem in $\R^n\setminus{\mathcal T}$,  whose solution, when unique,  simultaneously gives both  $\mathcal{V}$  and $Dom(\mathcal{V})$.  More precisely, let
\begin{equation}
\label{HJB2}\begin{array}{l}
{K}(x,u,p)
\doteq\max_{(w_0,w)\in S(A)}\{-\langle p,\overline f(x,w_0,w)\rangle-\overline l(x,w_0,w)+\overline  l(x,w_0,w)u \}. 
\end{array}
\end{equation}

%
 
\begin{Theorem}\label {ET3} Under the same hypotheses of Theorem \ref{EP3},
 there is a unique nonnegative viscosity solution ${\mathcal U}$  to 
\begin{equation}\label{BVPK}
\left\{\begin{array}{l}
{K}(x,u(x),Du(x))= 0\quad \qquad \text{in } \  \R^n\setminus\mathcal{T}\\
 u(x)= 0\qquad\qquad\quad\qquad \text{on } \ \partial\mathcal{T}.
\end{array}\right.
\end{equation}
 Moreover,  ${\mathcal V}\equiv V^r\equiv \Sigma\equiv\Psi^{-1}({\mathcal U})=-\log(1-{\mathcal U})$ and $Dom(\mathcal{V})=\{x:\ \ {\mathcal U}(x)<1\}$.
 
\noindent If we drop  {\rm (H0)$_1$},  {\rm (H0)$_2$} in (i), ${\mathcal U}$ (possibly $\ne \Psi(V^r)$) is the unique solution just  among the continuous functions. 
\end{Theorem} 
{\it Proof.}\, Let us prove the theorem in case  {\rm (H0)$_1$},  {\rm (H0)$_2$} are not assumed.  In order to apply the uniqueness result proved in Theorem  4.7 in \cite{MS}, let us observe that, under hypotheses   {\rm (H2)} and  {\rm (H3)},  the {\it asymptotic}  and the {\it minimal exit-time value functions}  ${\mathcal V}$ and ${\mathcal V}^m$,  as well as  their extended versions $V$ and  $V^m$ there  introduced,  do all coincide.  They also are equal to our infinite horizon value function ${\mathcal V}$ ($\equiv V$ by Theorem \ref{extT}).
Indeed, owing to  {\rm (H2)} and  {\rm (H3)}, both original and extended nearly optimal trajectories have to  approach at least asymptotically $\C$.
In fact,  since $V\equiv {\mathcal V}$, the conditions in hypothesis (H2) hold for $V$ too, and as shown  in the proof of Proposition \ref{Pu}, the liminf  in (\ref{H3}) is zero also for the $\varepsilon$-optimal trajectories of the original system.   Thanks to (\ref{CC}),  the last statement  follows now from  (i) of Theorem 4.7 in \cite{MS}, while the first statement  is a consequence of  (ii) of Theorem 4.7 in \cite{MS} together with either Theorem \ref{COVr} when (ii) is assumed or Theorem \ref{Thril}, when (i) holds.   \qed

 \begin{Remark}\label{BCs} {\rm 
Since when (H2) and (H3) hold for the same $\C$, the infinite horizon value function ${\mathcal V}$ coincides with the asymptotic exit-time value function considered in \cite{MS},   sufficient conditions for its continuity can be found there (see (TPK)$'$ in \cite{MS}).  In particular, when (H2) holds for $\C$, in view of Proposition 6.2 in \cite{MS}, (SC1) or (SC2) for the same $\C$ imply not only (H3), but also the  continuity of ${\mathcal V}$ and the boundary condition (\ref{Bcn}).  Moreover, as already observed, they also yield (H0)$_2$ (i).   Since in this section we suppose $l$ locally Lipschitz continuous in $x$, condition  (H0)$_2$ (ii) is trivially verified.}
\end{Remark} 

Therefore we have
 \begin{Corollary}\label{VinftyUn}  Let $\C\times\{0\}$ be a viability set for $(f,l)$.  Assume  the existence of a local MRF  and either {\rm (SC1)} or {\rm (SC2)} for ${\mathcal T}$. Then
 \begin{itemize}
 \item[(i)]
there is a unique nonnegative viscosity solution ${\mathcal U}$ to {\rm(\ref{BVPK})}, which turns out to be
continuous. Moreover,   ${\mathcal V}\equiv V^r\equiv\Sigma \equiv\Psi^{-1}({\mathcal U}) = -\log(1-{\mathcal U})$ and $Dom(\mathcal{V}) = \{x : {\mathcal U}(x) < 1\};$
 \item[(ii)] $\mathcal{V}$  $(\equiv V^r\equiv \Sigma)$   is the unique nonnegative  viscosity solution to  { \rm   (\ref{Ei})} in $Dom({\mathcal V})$  among the pairs $(u,\Omega)$ in ${\mathcal S}$.  Moreover,  $\mathcal{V}$ is continuous.
 \end{itemize}
\end{Corollary}

When $A$ is unbounded, the case  $q=p$ is  the only one in which we could have ${\mathcal V}(x)> V(x)$ for some $x$.  Since $\Sigma(x)=V^r(x)$, in order to characterize $\Sigma$, the well-posedness, that is the equality    ${\mathcal V}\equiv V$,  is not required. Hence in this  whole section assumption (H2) could be weakened,  by replacing in it the function  ${\mathcal V}$ with $V$. Accordingly, in Corollary \ref{VinftyUn}  it would be enough to assume $\C\times\{0\}$   viable   for $(\ol{f},\ol{l})$ and the existence of a MRF for the extended setting.
 
  
 \section{ Discounted infinite horizon approximations}\label{Dis}
 
In this section we give a representation formula for the limit as $\delta$ tends to $0^+$ of the infinite horizon value function with  discount rate $\delta>0$:
 $$
{\mathcal V_\delta}( x)\doteq\inf_{\alpha\in{\mathcal A} }\int_0^{+\infty}e^{-\delta\,t}l(y(\t),\alpha(\t))\,d\t.
$$ 
To this aim, for any $\delta>0$,  when $A$ is unbounded, we also introduce the extended value function
 $$
 V_\delta( x)\doteq\inf_{(w_0,w)\in\Gamma}\int_0^{+\infty}e^{-\delta\int_0^sw_0^q(s)\,ds}\ol{l}(\xi(s), w_0(s),w(s))\,ds,
$$ 
and,  agreeing with the notation of Subsection \ref{Rel},  if $A$  is compact [resp., unbounded],  we consider  the relaxed version of ${\mathcal V}_\delta$, ${\mathcal V}^r_\delta$ [resp., of $V_\delta$,  $V_\delta^r$].
 
As a first step, by Proposition 3.2 in \cite{M}  all these value functions  are  supersolutions to
\bel{delta}
\delta u+{\mathcal H}(x, Du(x))=0 
\eeq
in $\R^n$.  If they are locally bounded and with open domains, they also are subsolutions to (\ref{delta}) in their domains.   
Notice that, when $A$ is unbounded, by Theorem 2.1 in \cite{M}, equation (\ref{delta}) can be replaced by 
$$
H_\delta(x,u(x), Du(x))=0 \qquad x\in\R^n,
$$
where, for any $(x,r,p)\in\R^{2n+1}$,  $H_\delta$ is the following  continuous Hamiltonian 
\bel{deltae}
H_\delta(x,r,p)\doteq \max_{(w_0,w)\in S(A)
}\left\{\delta r\,w_0^q-\langle\overline{f}(x,w_0,w),\, p\rangle-\overline{l}(x,w_0,w)\right\}.
\eeq
By Corollary 4 in \cite{MS2},  for any $\delta>0$ we have what follows. 
 \begin{Theorem}\label{Udelta}  If ${\mathcal V}_\delta$ is bounded, then it is the unique bounded  solution to  {\rm (\ref{delta}) } in $\R^n$ and it is continuous.  Hence, if $A$ is compact one has  ${\mathcal V}_\delta\equiv {\mathcal V}^r_\delta,$ and   ${\mathcal V}_\delta\equiv V_\delta\equiv   V^r_\delta$ otherwise.
 \end{Theorem}

\begin{Remark}\label{csdelta}{\rm  It is easy to see that, when $A$ is unbounded,  sufficient conditions in order to have ${\mathcal V}_\delta$ bounded are, for instance, either  
$$
\begin{array}{l}
|f(x,a)|\le \bar M +M(1+|x|)|a|^p \quad \text{and} \quad l(x,a)\le \bar  M(1+|x|^r)+ M_R|a|^q \\
\text{or } \\
l(x,a)\le \bar M + M_R|a|^q   \quad\text{ $\forall (x,a)\in\R^n\times A$ with  $|x|\le R$,}
\end{array}
$$
for some $\bar M>0$, $r\ge1$  ($M_R$ is the same as in (\ref{hp0})).  Formally, the same conditions  with $a=0$  yield the boundedness of ${\mathcal V}_\delta$ for $A$ bounded. 

\noindent We refer to Corollary 4 in \cite{MS2}, for a  characterization of  ${\mathcal V}_\delta$  as unique solution  to (\ref{delta}) in $\R^n$   in some classes of unbounded functions with prescribed growth at infinity. 
}
\end{Remark}

\begin{Theorem}\label{supdelta}  Assume that each  ${\mathcal V}_\delta$ is bounded. Then
$$
\lim_{\delta\to0^+}{\mathcal V}_\delta(x) =V^r(x) \qquad \forall x\in\R^n.
$$
\end{Theorem}

\noindent{\it Proof.}  We give the proof in the case $A$ unbounded, being the other case similar. Taking into account that  the sequence $\delta\to {\mathcal V}_\delta$ is monotone non increasing, by Theorem \ref{Udelta}, we have
$$
\Lambda(x)\doteq \lim_{\delta\to0^+}{\mathcal V}_\delta(x)=\sup_{\delta>0} {\mathcal V}_\delta(x)=\sup_{\delta>0} V^r_\delta(x)\le V^r(x)
$$
for every $x\in\R^n$.  In view of Theorem \ref{5.2} (ii),  $V^r$ is the minimal supersolution to  (\ref{Ei}) in $\R^n$, hence  it is now sufficient to show that $\Lambda$ ($=\Lambda_*$) is a supersolution to (\ref{Ei}) in $\R^n$ for any $x$ such that    $\Lambda(x)<+\infty$.

\noindent By the monotonicity of the sequence  ${\mathcal V}_\delta$ and by the continuity of  each   ${\mathcal V}_\delta$, it is known that   $\Lambda(x)=\Lambda_*(x)=\underset{\delta\to0^+}{\lim\inf_*}{\mathcal V}_\delta(x)$ (see \cite{BCD}). The claim  follows now from stability results of viscosity solutions, taking into account the continuity of the  ${\mathcal V}_\delta$ and the fact that we can consider  the regular Hamiltonian in (\ref{deltae}). \qed
  
  \vv
  In the above proof we used the upper optimality principle. Of course, it is also possible to obtain it by working directly on the control problem.

\section{Ergodic problem}\label{erg}
In this section we briefly investigate the so-called ergodic problem, that is  the convergence of the limits $\underset{t\to+\infty}\lim{\mathcal V} (t,x)/t$, $\underset{\delta\to0^+}\lim\delta\,{\mathcal V}_\delta(x)$.     Our goal here is just to describe  how known hypotheses and proofs can be adapted to the case of unbounded controls. Hence in the sequel we consider $A$ unbounded and  assume 
$f$ and $l$ periodic in the state variable and  global controllability. 
Our precise assumptions, together with (H1),  are the following.  
{\it \begin{itemize}
\item[{\bf (H4)}]  (i)\,  $T_i>0$ \,$(i=1,\dots,n)$  are real numbers and the functions $f(x,a)$, $l(x,a)$ are periodic in $x_i$ with the period $T_i$ \,$(i=1,\dots,n)$. Moreover  there are $L$ and $M>0$ such that $\forall x$, $x_1$,  $x_2\in{\bb T}^n$, $\forall a\in A$,
\begin{equation}\label{hp03}
\aligned
|f(x_1,a)-f(x_2,a)|&\le L (1+|a|^p)|x_1-x_2|,\\
 |l(x_1,a)-l(x_2,a)|&\le L(1+|a|^q)|x_1-x_2| \\
  l(x,a) \le M(1+|a|^q),& \quad
 |f(x,a)|\le M(1+|a|^p),
\endaligned
\end{equation}
 where
  ${\bb T}^n$ denotes the  $n$--dimensional torus  $\R^n\,/\, (\Pi_{i=1}^nT_i{\Bbb Z})\,\sim\, \Pi_{i=1}^n[0,T_i]$.

  (ii) There are $C$, $\gamma>0$ such that for any pair $x$, $z\in {\bb T}^n$ there exist 
$S>0$ and   $\mu\in\Gamma^r$ such that $\xi^r_x(S, \mu)=z$ and $S\le C|x-z|^\gamma$.
\end{itemize}}

A sufficient condition to have (H4) (ii) (with  $\gamma=1$) is the usual hypothesis that, for some $r>0$,  
$B(0,r)\subset\overline{\text{co}}\, \ol{f}(x, S(A))$ for any $x\in\R^n.$

%
%

\begin{Remark}\label{Noerg}{\rm  Owing to Theorems \ref{supC} and \ref{supdelta}, at least when any ${\mathcal V}_\delta$ is bounded,  $\underset{t\to+\infty}\lim{\mathcal V} (t,x)= \underset{\delta\to0^+}\lim {\mathcal V}_\delta(x)=V^r(x)$ for every $x\in\R^n$.  As a consequence,  the  limits  $\underset{t\to+\infty}\lim{\mathcal V} (t,x)/t$ and $\underset{\delta\to0^+}\lim\delta\,{\mathcal V}_\delta(x)$ converge  obviously  to  zero  when  $V^r$ is finite in $\R^n$.  In fact, being  $l\ge0$  such a convergence is  locally uniform. 

\noindent When $l\le M(1+|a|^q)$ and   (H4) (ii) is in force,   $V^r$ is finite as soon as   $(f ,l)(x,a)=(0,0)$ for some pair $(x,a)$, or, more in general, if there exists a subset $\C\subset{\mathcal Z}$ such that $\C\times \{0\}$ is a viability set for $(f,l)$. In this case indeed, for any $x\in\R^n$ it is possible to construct an admissible control $\alpha$ with  finite cost,  by concatenating a control steering $x$ to $\C$ in time $T$, as in (H4) (ii),    with a control keeping the trajectory inside $\C$ with null cost for all  $t>T$. Such a control  exists in view of the viability assumption.
%
}
\end{Remark}

\begin{Proposition}\label{Vdelta}
Assume {\rm(H3)}. Then, for any $x$, $z\in {\bb T}^n$,
\bel{esterg}
0\le {\mathcal V}_\delta(x)\le M/\delta, \qquad
|{\mathcal V}_\delta(x)-{\mathcal V}_\delta(z)|\le M C\,|x-z|^\gamma.
\eeq
 Moreover, setting
${\mathcal W}_\delta(x)\doteq {\mathcal V}_\delta(x)-{\mathcal V}_\delta(0)$, one also has
\bel{estergW}
|{\mathcal W}_\delta(x)|  \le M_1, \qquad
|{\mathcal W}_\delta(x)-{\mathcal W}_\delta(z)|\le M C\,|x-z|^\gamma,
\eeq
where $M_1\doteq MC(\sqrt{n}\max_{i=1,\dots,n}T_i)^\gamma$.
%
\end{Proposition}
{\it Proof.} In view of Theorem \ref{Udelta}, for any $x\in\R^n$  one has ${\mathcal V}_\delta(x)= V_\delta(x)\equiv V^r_\delta(x)$. Therefore the first estimate in (\ref{esterg}) follows immediately from the fact that $\ol{l}\le M$, considering the relaxed control $\mu\equiv0$.  Assuming $V^r_\delta(x)-V^r_\delta(z)\ge 0$, as it is not restrictive,  the second inequality in  (\ref{esterg}) can be obtained  plugging in the DDP for $V^r_\delta(x)$ the control given by (H4)  (see e.g. Theorem 2 in \cite{A}). Both the estimates in (\ref{estergW}) are easy consequence of  (\ref{esterg}). \qed

\begin{Theorem}\label{Therg} Assume {\rm (H4)}. Then there exists a constant  $\lambda\ge 0$ such that 
$$
 \underset{\delta\to0^+}\lim\delta\,{\mathcal V}_\delta(x)=\lambda, \quad \underset{t\to+\infty}\lim{\mathcal V} (t,x)/t=\lambda \quad \text{ uniformly in $\R^n$.}
$$
  Moreover, there exists some  $\delta_n\to0^+$ such that
$$
 \underset{n\to+\infty}\lim {\mathcal W}_{\delta_n}(x)={\mathcal W}_0 \quad\text{ uniformly in $\R^n$,} 
$$
  and ${\mathcal W}_0\in BUC(\R^n)$ is a solution of
 \bel{Herg0}
 \tilde H_\lambda(x, Du(x))=0 \qquad \text{in $\R^n$,}
\eeq 
where
$$
\tilde H_\lambda(x, p) \doteq\max_{(w_0,w)\in S(A)}\left\{-\langle\overline{f}(x,w_0,w),\,p\rangle-\overline{l}(x,w_0,w)+\lambda w_0^q\right\}.
$$
\end{Theorem}

\noindent{\it Proof.}
By Proposition \ref{Vdelta}, the Ascoli-Arzel\`a Theorem and  the periodicity of the solutions imply that  there exists a sequence $\delta_n\to0^+$ such that $ \underset{n\to+\infty}\lim\delta_n{\mathcal V}_{\delta_n}=\lambda\in{\mathcal C}(\R^n)$ and $\underset{n\to+\infty}\lim {\mathcal W}_{\delta_n}={\mathcal W}_0\in{\mathcal C}(\R^n)$. The second inequality in
(\ref{esterg}) implies that $\lambda$ is a constant and consequently $\delta_n{\mathcal W}_{\delta_n}\to 0$ uniformly in 
$\R^n.$ It is now easy  to check that ${\mathcal W}_{\delta}$ satisfies
$$\max_{(w_0,w)\in S(A)
}\left\{\delta {u}\,w_0^q-\langle\overline{f}(x,w_0,w),\,D{u}\rangle-\overline{l}(x,w_0,w)+\delta{\mathcal V}_{\delta}(0)w_0^q\right\}=0.$$
By the stability of the viscosity solutions  and by the regularity of the above  Hamiltonian,  
it follows that  $(\lambda,{\mathcal W}_0)$ solves $\tilde H_\lambda(x, Du)=0.$
It remains to be proved that $\lambda$ is uniquely determined and that the whole family 
$\delta{\mathcal V}_\delta$ converges to $\lambda.$ The claim is that there exists a unique $\lambda\ge 0$ such that
(\ref{Herg0}) has a bounded, uniformly continuous solution in $\R^n$. 
First let us prove that if there exist $\lambda_1,$ $\lambda_2\ge 0$, such that  $u_1$ is a subsolution
to $\tilde H_{\lambda_1}(x, Du) =0$ and $u_2$ is a supersolution to  $\tilde H_{\lambda_2}(x, Du) =0$ then one must have
$\lambda_1\le\lambda_2.$ Let us argue by contradiction and assume $\lambda_1>\lambda_2.$
We can suppose, eventually adding a constant, that $u_1>u_2$. Let $\varepsilon$ be small enough such that
$\lambda_1-\varepsilon u_1>\lambda_2-\varepsilon u_2$ in $\R^n.$ Therefore $u_2$ is also a supersolution to 
$$ \max_{(w_0,w)\in S(A)}
\left\{\varepsilon u_2w_0^q-\langle\overline{f}(x,w_0,w),\,Du_2\rangle-\overline{l}(x,w_0,w)+
(\lambda_1-\varepsilon u_1)w_0^q\right\}=0$$
and
$u_1$ is also a subsolution to 
$$ \max_{(w_0,w)\in S(A)}
\left\{\varepsilon u_1w_0^q-\langle\overline{f}(x,w_0,w),\,Du_1\rangle-\overline{l}(x,w_0,w) 
+(\lambda_1-\varepsilon u_1)w_0^q\right\}=0$$
in $\R^n$. By the comparison principle underlying Theorem \ref{Udelta}  we would get $u_1(x)\le u_2(x),$ a contradiction.  Therefore the claim is proved and  one has $\lambda_1\le\lambda_2$.

\noindent  Now let us assume that there exist $\lambda_1=\lim_{\delta_n\to 0}{\mathcal V}_{\delta_n}$ and
 $\lambda_2=\lim_{\bar \delta_n\to 0}{\mathcal V}_{\bar \delta_n}$. The above result yields that $\lambda_1=\lambda_2$, so that the uniform  limit $ \underset{\delta\to0^+}\lim\delta\,{\mathcal V}_\delta(x)=\lambda$ is proved. 
 
 In order to prove that $ \underset{t\to+\infty}\lim{\mathcal V} (t,x)/t=\lambda$ uniformly, for the same $\lambda$ as above, let us first introduce the function 
 $ v(t,x)\doteq C+  {\cal W}_0(x)+\lambda t$ for all $(t,x)\in\R_+\times\R^n, $
where ${\cal W}_0$ is a solution to  $\tilde H_\lambda(x, Du)=0$ and $C>0$ is chosen so that $C+{\cal W}_0\ge 0$. Then $v$ is a supersolution to (\ref{CPb}) for any $T>0$ and by the comparison principle underlying Theorem \ref{P1b}, 
$$
 {\mathcal V} (t,x)\le v(t,x)=  C+  {\cal W}_0(x)+\lambda t \qquad \forall (t,x)\in\R_+\times\R^n.
$$
 Let us now consider the function $\tilde v(t,x)\doteq  -C+  {\cal W}_0(x)+\lambda t$ for all $(t,x)\in\R_+\times\R^n, $
where and $-C+{\cal W}_0\le 0$. Then $v$ is a subsolution to (\ref{CPb}) for any $T>0$ and  we get
$$
 {\mathcal V} (t,x)\ge \tilde v(t,x)= - C+  {\cal W}_0(x)+\lambda t \qquad \forall (t,x)\in\R_+\times\R^n,
$$
arguing as above.  By the last two inequalities, the proof follows. \qed

\begin{Remark}\label{Vnb} {\rm Let us observe that  the effective Hamiltonian $\tilde H_\lambda$ really determines $\lambda$.  This would not be the case, if there existed a function ${\cal W}_0\in BUC(\R^n)$ such that   the max in  the definition of $\tilde H_\lambda$ was reached for every $x\in\R^n$ in a vector $(0,w)\in S(A).$ If fact, such a function would be a solution of  
$$
\max_{(0,w)\in S(A)}\left\{-\langle\overline{f}(x,0,w),\,Du\rangle-\overline{l}(x,0,\bar w)\right\}=0,
$$ 
 and then it would also solve $\tilde H_\lambda(x, Du)=0$ for all $\lambda$. However, applying Theorem \ref{5.2}, such ${\cal W}_0$ would be greater than 
  the value function of an infinite horizon problem with compact controls $(0,w)\in S(A)$ (where $|w|^q=1$) and  lagrangian $\bar l(x,0,\bar w)\ge C_2$, equal to $+\infty$.  Again, the coercivity hypothesis (\ref{CC}) plays   a crucial role.}
 \end{Remark}

 \end{large}
\end{document}